\documentstyle[11pt]{article}

\newlength{\abstractwidth}
\renewcommand{\thefootnote}{\fnsymbol{footnote}}
\renewcommand{\thanks}[1]{\footnote{#1}} 
\newcommand{\starttext}{ \setcounter{footnote}{0}
\renewcommand{\thefootnote}{\arabic{footnote}}}

\newcommand{\be}{\begin{equation}}
\newcommand{\bea}{\begin{eqnarray}}
\newcommand{\eea}{\end{eqnarray}} \newcommand{\ee}{\end{equation}}
 \newcommand{\<}{\langle}
\renewcommand{\>}{\rangle} \def\ba{\begin{eqnarray}}
\def\ea{\end{eqnarray}}

\def\D{{\cal D}}

\def\cL{{\cal L}}

\def\r{\rho}

\def\tr{{\rm tr}}
\def\det{{\rm det}}

\def\half{ {1\over 2}}

\def\log{\,{\rm log}\,}
\def\exp{\,{\rm exp}\,}

\def\a{\alpha}

\def\g{\gamma}
\def\d{\delta}
\def\e{\varepsilon}

\def\l{\lambda}

\def\f{\phi}
\def\r{\rho}
\def\si{\sigma}

\def\D{\Delta}
\def\O{\Omega}

\def\na{\nabla}

\def\ve{\varepsilon}
\def\ge{\geq}
\def\le{\leq}
\def\ddt{{\partial\over\partial t}}

\def\ti{\tilde}

\def\R{{\bf R}}

\def\i{\infty}

\def\p{\partial}
\def\s{\sum}

\def\D{\Delta}

\def\cF{{\cal F}}

\def\cW{{\cal W}}

\def\na{{\nabla}}

\def\[{{\bf [}}
\def\]{{\bf ]}}

\newcommand{\abs}[1]{|#1|^2}

\def\XXint#1#2#3{{\setbox0=\hbox{$#1{#2#3}{\int}$ }
		\vcenter{\hbox{$#2#3$ }}\kern-.6\wd0}}

\newtheorem{prop}{Proposition}

\newcommand{\ric}{{Ric}}
\newcommand{\innpro}[1]{\langle#1\rangle}
\newcommand{\bk}[1]{\Big(#1\Big)}

\makeindex

\newcommand{\sic}{{Sic}}

\begin{document}
\starttext \baselineskip=18pt \setcounter{footnote}{0}
\newtheorem{theorem}{Theorem}
\newtheorem{lemma}{Lemma}
\newtheorem{corollary}{Corollary}
\newtheorem{definition}{Definition}
\newtheorem{conjecture}{Conjecture}
\newtheorem{proposition}{Proposition}

\begin{center}
{\Large \bf PSEUDO-LOCALITY FOR A COUPLED RICCI FLOW
\footnote{   AMS classification numbers: 53C23, 53C44. Work supported in part by NSF grant DMS-12-66033.}}

\bigskip

{\large Bin Guo, Zhijie Huang, Duong H. Phong } \\

\medskip

\begin{abstract}

\medskip
\small{
Let $(M,g,\f)$ be a solution to the Ricci flow coupled with the heat equation for a scalar field $\f$. We show that a complete, $\kappa$-noncollapsed solution $(M,g,\f)$ to this coupled Ricci flow with a Type I singularity at time $T<\i$ will converge to a non-trivial Ricci soliton after parabolic rescaling, if the base point is  Type I singular. A key ingredient is a version of Perelman pseudo-locality for the coupled Ricci flow.
}

\end{abstract}

\end{center}

\baselineskip=15pt
\setcounter{equation}{0}
\setcounter{footnote}{0}

\section{Introduction}
\setcounter{equation}{0}

The Ricci flow \cite{H} can be viewed as the parabolic and Euclidean version of Einstein's equation in the vacuum. In presence of matter fields, Einstein's equation becomes a coupled system. Thus we should also consider the Ricci flow coupled with other flows. The simplest is the Ricci flow coupled with the heat equation for a scalar field. 
This is a special case of the Ricci flow coupled with the harmonic map flow, a version of which had been instrumental in the proof of the short-time existence of solutions to the Ricci flow \cite{DeT}. Coupled Ricci flows also arise as dimensional reductions of the Ricci flow in higher dimensions, and the coupling to a scalar field arises in particular in the Ricci flow on warped products \cite{LS}.

\medskip
More specifically, let $M$ be a compact manifold. The Ricci flow coupled to the heat equation is the following system of equations for a metric $g_{ij}(t)$ and scalar field $\f(t)$,
\be \label{RFHE}
\left\{
\begin{array}{l}
	\frac{\p g}{\p t} = -2Ric_g + 2d\f \otimes d\f, \qquad 
	\f_t = \Delta_g \f \\
	g(0) = g_0, \qquad \f(0)=\f_0
\end{array}
\right.
\ee
where $g_0$ and $\f_0$ are given smooth initial data, and
the coupling constant to the scalar field has been normalized to be $1$. Its stationary points are solutions to Einstein's equation
with $\f$ the matter field and stress tensor $T_{ij}=\p_i\f\p_j\f$.
This flow has been first studied extensively by List \cite{L}, who established criteria for its long-time existence, and obtained extensions to this case of Perelman's monotonicity formula and non-collapse results, as well as an extension of Hamilton's compactness theorem. Similar results for the more general case of the Ricci flow coupled with the harmonic map flow were subsequently obtained by M\"uller \cite{M}.

\medskip
The goal of this paper is to establish a Perelman pseudo-locality theorem for the Ricci flow coupled with a scalar field. Let
\bea
Sic_{g,\f} = Ric_g -d\f\otimes d\f
\eea
and denote its components by $S_{ij}=R_{ij}-\p_i\f\p_j\f$
and its trace by $S =\tr_g Sic_{g,\f} = g^{ij} S_{ij} = R_g-|\na \f|^2_g$.
From now on, we shall omit the sub-script $g$ and $\f,$ if it is clear from the context.
We prove the following pseudo-locality theorem:

\begin{theorem}
\label{1} {\rm [Pseudo-locality]} Given $\a \in (0, \frac{1}{100 n})$, there exist $\ve =\ve(n,\a,C), \d = \d(n,\a,C)$ with the following property.
For any solution to the Ricci flow coupled with a scalar field $(M,g(t),\f(t), p)$,
$t\in [0,(\ve r_0)^2]$, which has complete time slices and satisfies

\smallskip
{\rm (1)} $S(g(0))\geq - r_0^2$ on the ball $B_{g(0)}(p,r_0)$; 

{\rm (2)} $Area_{g(0)}(\p \O)^n \geq (1-\d)c_n Vol(\O)^{n-1}$, for any $\O \subset B_{g(0)}(p,r_0)$, where $c_n$ is the isoperimetric constant in $\R^n$,

{\rm (3)} $|\f_0| \leq C$ on the ball $B_{g(0)}(p,r_0)$ for some constant $C$,

\smallskip
\noindent
we have 
	\be
	|Rm|(x,t) \leq \alpha t^{-1} +(\ve r_0)^{-2} ,
	\ee
	for any $(x,t)$ such that $d_{g(t)}(x,p) \leq \ve r_0$ and $t\in ( 0 , (\ve r_0)^2].$
\end{theorem}

A well-known conjecture of Hamilton is that blow-ups of Type I singularities in the Ricci flow should converge to a non-trivial gradient soliton. This conjecture was proved in dimension $n=3$ by Perelman \cite{P}. The convergence to a gradient soliton for a Type I singularity was proved in all dimensions by Naber \cite{N}, and the non-triviality of the soliton as $t$ tends to the maximum existence time was subsequently proved by Enders, M\"uller, and Topping \cite{EMT}. 
In a different direction,
it was also shown by Cao and Zhang \cite{CZ}
that the blow-down limit of Type I $\kappa$-noncollapsed ancient solutions
was a non-trivial soliton.
A key ingredient in 
the arguments of Enders, M\"uller, and Topping was Perelman's pseudo-locality theorem. In this paper, we extend their arguments to the coupled Ricci flow as follows. First, List \cite{L} has shown that the maximum existence time $T$ for the coupled Ricci flow must satisfy
\bea
{\rm lim}_{t\to T}\,{\rm sup}_{x\in M}|Rm|^2(t,x)=\infty.
\eea
Thus, in analogy with the Ricci flow, a solution $(M,g, \f)$ of the coupled Ricci flow with maximal time $T$ 
is called of Type I if there exists some constant $C_0$ such that 
\be
\label{Type I} 
	\sup_{x\in M} |Rm_{g(t)}|(x,t) \leq \frac{C_0}{T-t}, \quad t\in [0,T).
\ee 
A point $p$ is said to be a Type I singular point if there exists a sequence $(p_i,t_i)$, $p_i \to p, t_i\to T$, with $|Rm_{g(t_i)}|(p_i)\geq c(T-t_i)^{-1}$ for some constant $c>0$.
Next, the notion of gradient soliton for the coupled flow can be extended as a triple $(g_{ij},\f, f)$ satisfying
\bea
Sic+\na^2 f-{g\over 2(T-t)}=0, \quad \Delta \phi-\<\na f,\na \phi\>=0.
\eea
The soliton is said to be trivial if the metric $g_{ij}$ is flat.
As a consequence of the above pseudo-locality theorem, we 
have then

\begin{theorem}
\label{2}
	Let $(M, g(t),\f(t))$ be a solution of the coupled Ricci flow (\ref{RFHE}) with $|\phi_0|\leq C$, and assume that it has a Type I singularity at time $T<\i$, with $p$ a Type I singularity point. Let
$\l_i \to \i$ be any sequence of numbers, and define a sequence of coupled Ricci flows by 
	\be 
		g_i(t) = \l_i g(\l_i^{-1} t + T), \ \ \f_i(t) = \f (\l_i^{-1} t + T), \ \ \forall t\in [-\l_i T,0)
	\ee
	Then there exists a subsequence of $(M,g_i,\f_i,p)$ which converges to a non-trivial gradient shrinking Ricci soliton $(M_\i,g_\i(t), \f_\i, p_\i)$. The function $\f_\i$ is actually constant, so that this soliton for the coupled flow actually reduces to a soliton for the usual Ricci flow.
\end{theorem}

\section{The evolution equation for the curvature}
\setcounter{equation}{0}

For the convenience of the reader, we quote here several formulas and estimates already established by List \cite{L} and M\"uller \cite{M}.

\smallskip
Let $(M,g_{ij}(t),\f(t))$, $t\in [0,T)$, be a solution of the coupled Ricci flow (\ref{RFHE}). Then
\cite{L}
\be 
	\ddt S_{jk} = \D S_{jk}+ 2 R_{ijk}{}^l S_{l}{}^i - S_{j}{}^l R_{lk}- S_k{}^l R_{lj} + 2 \D \f \cdot \f_{jk}.
\ee 
Taking traces with respect to $g^{jk}$ gives
\be 
	\ddt S = \D S + 2 |S_{ij}|^2 + 2 (\D \f)^2 
\ee

List \cite{L} showed $\f(x,t)$ is uniformly bounded along the flow
\begin{lemma}(Lemma 5.10 in \cite{L}) Let $(g(t),\f(t))$ be a solution to the coupled Ricci flow (\ref{RFHE}) on $M\times [0,T)$ with initial data $(g_0,\f_0)$, assume $\sup|\f_0|\leq C$, then we have for $t>0$, 
	\bea
		\inf_{x \in M}\f_0(x) \leq \f(x,t) \leq \sup_{x\in M} \f_0(x) \label{List estiamte C0}\\
		\sup_{x\in M} |\na \f|^2(x,t) \leq C^2t^{-1},  \label{List estimate 0}
	\eea 	
\end{lemma}

The following derivative estimate hold
\begin{prop}(Theorem 5.12 in \cite{L})
Let $(g(t),\f(t))$ be a solution to the coupled Ricci flow (\ref{RFHE}). Fix $x_0\in M$ and $r>0$, if
\be 
\sup_{B(T,x_0,r)} r^2 |Rm| \leq \ti C
\ee 
where  $B(T,x_0,r)$ is the geodesic ball centered at $x_0\in M$ with radius $r$ at time $T$.  Denote $\Phi = (Rm, \na^2 \f)$,
then the derivatives of $\Phi$ satisfy for all $m\geq 0$ and for all $t\in(0,T]$ the estimates hold
\be 
\label{List derivative estimate}
\sup_{B(t,x_0,r/2)}  |\na^m \Phi|^2 \leq C(n,m) \ti C^{m+2} \left( r^{-2} +t^{-1} \right)^{m+2}
\ee 
where $C=C(n,m)$ is a constant depending only on $n$ and $m$. 

\end{prop}

\section{Reduced distance and volume}
\setcounter{equation}{0}

In this section we provide some background material on the reduced distance and volume for the coupled Ricci system. This material can be found in M\"uller \cite{M} and Vulcanov \cite{V}, but we provide a self-contained discussion for completeness and for easier use when it is needed in Sections 4-6.

\medskip
 
Let $(g_{ij}(t),\f(t))$ be a solution to the coupled Ricci flow (\ref{RFHE}) on $t\in (0,T)$. For some fixed $t_0 \in [0,T)$, set $\tau = t_0 - t$. In terms of $\tau$, the flow becomes 
\bea 
(g_{ij})_\tau =2S_{ij},
\qquad
\f_\tau =-\D \f.
\eea

\subsection{The $\cL$-length of a path for coupled Ricci flow}

Let $\g: [\tau_1,\bar\tau] \to M$ be a path, where $0\leq \tau_1$. We define the $\cL$-length of the path $\gamma$ by 
\be
	\cL(\g) = \int_{\tau_1}^{\bar\tau}\sqrt{\tau} \Big( S(\g(\tau)) + |\g'(\tau)|^2 \Big) d \tau 
\ee
where $S(\g(\tau))$ and the norm $|\g'(\tau)|$ are evaluated using the metric $g_{ij}(t)$ at time $t = t_0 - \tau$.

\subsubsection{The first variational formula for $\cL$}

Let $\g_u(\tau)$ be a family of paths in $M$ with tangent vector field $X$ and variational vector field $Y$, i.e.,
\bea
X(\tau)=\frac{\p\g_u}{\p\tau}(\tau),
\qquad
Y(\tau) = \frac{\p \g_u}{\p u}(u,\tau).
\eea
Note that $\na_Y X = \na_X Y$. We assume momentarily that $\tau_1>0$.
Now
\bea 
\label{variationcL1}
   \d_Y \cL(\g)	&=& \frac{d}{du} \Big|_{u=0} \int_{\tau_1}^{\bar\tau}\sqrt{\tau} \Big( S(\g_u(\tau)) + |\g_u'(\tau)|^2 \Big) d \tau \nonumber \\
				&=&  \int_{\tau_1}^{\bar\tau} \sqrt{\tau} \Big( \< \na S(\tau) , Y(\tau)\> + 2 \< \na_Y X, X\> \Big) d \tau.  
\eea
Since $\frac{d}{d\tau}\< Y, X \> =2Sic(Y,X)+ \<\na_X Y, X \>+ \< Y, \na_X X\>$, the second term on the right hand side can be rewritten as
\bea 
\int_{\tau_1}^{\bar\tau} \sqrt{\tau} \<\na_X Y, X \> d\tau
	&=& \int_{\tau_1}^{\bar\tau} \sqrt{\tau} \Big(\frac{d}{d\tau}\< Y,X\> -\<\na_X X,Y\> - 2Sic(Y,X) \Big) d\tau \nonumber \\
	&=& \sqrt{\tau} \< Y, X\> \Big|_{\tau_1}^{\bar\tau} - \int_{\tau_1}^{\bar\tau} \sqrt{\tau} \<Y, \frac{1}{2\tau}X + \na_X X +2Sic(X,\cdot) \> d\tau. 
\eea 
Thus we obtain the following first variational formula
\be 
\d_Y \cL(\g) = 2\sqrt{\tau} \<X,Y\>\Big|_{\tau_1}^{\bar\tau} -  \int_{\tau_1}^{\bar\tau} \<Y, 2 \na_X X - \na S + \frac{1}{\tau}X + 4Sic(X,\cdot) \> d\tau 
\ee 

We define $\cL$-geodesics to be the critical points of $\cL$. The first variational formula gives immediately the equation for $\cL$-geodesics
\be 
\label{L-geod}
\na_X X - \frac{1}{2} \na S + \frac{1}{2\tau} X + 2 Sic(X, \cdot) = 0
\ee
where the connection and curvature are taken at the corresponding time $t= t_0 - \tau$ and the 1-form $Sic(X,\cdot)$ is identified with the vector field by the corresponding metric. 

\medskip
The $\cL$-geodesic can be rewritten as follows. 
Let $s= \sqrt{\tau}$ and $\ti \g(s) = \g(\tau(s)) = \g(s^2)$. Then setting $\bar s=\sqrt{\bar \tau}$ and $\ti X(s) =\ti \g'(s) = 2\sqrt{\tau} X(\tau)$, the geodesic equation becomes
\be
	\na_{\ti X} \ti X - 2s \na S + 4s^2 Sic(\ti X, \cdot) = 0
\ee
Henceforth, we assume that the paths $\g$ can be extended smoothly to $s=0$.
With the initial data $2v= \ti\g'(0)$, we can always solve for the geodesic equation for at least a short time, yielding a geodesic $\ti \g (s)$. Note also that the $\cL$-length becomes, in terms of the parameter $s$, 
\be
	L(\ti\g(\bar s), \bar s) = \int_0^{\bar s} \left(\frac{1}{2}|\ti \g'(s)|^2 +2 s^2 S(\ti \g(s)) \right)ds. 
\ee

For fixed $p\in M$, the $\cL$-exponential map is then defined as
the map $\cL \exp_{\tau} : T_p M \to M$ sending $v$ to $\g(\tau)$,
where $\g$ is the $\cL$-geodesic with $s=0$ and initial vector
$v = \lim\limits_{\tau \to 0} \sqrt{\tau} X(\tau) $.

\medskip
Next, the vector field $Y$ is said to be an $\cL$-Jacobi field if $\g_u(\tau)$ is an $\cL$-geodesic for all sufficiently small $u$. To get the equation for $\cL$-Jacobi fields, we assume that $X$ satisfies the $\cL$-geodesic equation for each $u$, and we differentiate with respect to $u$,
\be 
	\na_Y	 \big( 
				\na_X X - \frac{1}{2} \na S + \frac{1}{2\tau} X + 2 Sic(X, \cdot)  
			\big) = 0,
\ee
that is for any $W$, 
$$\< \na_Y	 \big( 
				\na_X X - \frac{1}{2} \na S + \frac{1}{2\tau} X + 2 Sic(X, \cdot)  
			\big), W \> =0.$$
Since 
\bea
&&
\na_Y \na_X X = \na_X \na_Y X -R(X,Y)X  = \na_X \na_X Y -R(X,Y)X
\nonumber\\
&&
\< \na_Y(\na S) , W \> = Y\< \na S, W \>  - \< \na S, \na_Y W\> = \na^2S(Y,W),
\nonumber\\
&&
\< \na_Y Sic(X, \cdot), W \> = Y \< Sic(X,\cdot),W \> - \< Sic(X,\cdot), \na_Y W\> = (\na_Y Sic)(X,W)+ Sic(\na_Y X, W)
\nonumber
\eea
we find that the $\cL$-Jacobi field equation can be written as
\be
\label{Jacobi}
T(Y)=0
\ee
where the vector field $T(Y)$ is defined by
\be 
\label{TY}
-T(Y)=  \na_X \na_X Y - R(X,Y)X +\frac{1}{2\tau }\na_Y X  -\frac{1}{2}\na^2S(Y,\cdot) + 2(\na_Y Sic)(X, \cdot)+ 2 Sic(\na_Y X, \cdot).
\ee 

\subsubsection{The second variational formula for $\cL$}

We now work out the second variation $Q(Y,Y)$ of the $\cL$-length, defined by
\bea
Q(Y,Y)=\delta_Y^2 \cL-\delta_{\na_YY}\cL.
\eea
Differentiating the formula (\ref{variationcL1}) for the first variation gives
\bea
	\d_Y^2 \cL 	&=& \int_0^{\bar \tau} \sqrt{\tau} \Big( 
					Y(Y S)+ 2\< \na_Y \na_Y X, X\> + 2 \< \na_Y X, \na_Y X\> \Big) d\tau  \nonumber \\
				&=&  \int_0^{\bar \tau} \sqrt{\tau} \Big( 
					Y(Y S) + 2 \< \na_Y \na_X Y ,X \> + 2 |\na_Y X|^2 \Big) d\tau  \nonumber \\
				&=& \int_0^{\bar \tau} \sqrt{\tau} \Big( 
					Y(Y S) + 2 \< \na_X \na_Y Y ,X \>+ 2 \<R(Y,X)Y,X \> + 2 |\na_X Y|^2 \Big) d\tau  
\eea 
Using the formula for Levi-Civita connection, we have 
\bea 
	&&\frac{d}{d\tau} \< \na_Y Y, X\> = \< \na_X \na_Y Y,X\> + \<\na_Y Y, \na_X X\> + 2Sic(\na_Y Y, X) + \< \dot{\na}_Y Y, X\> \nonumber \\
	&=& \< \na_X \na_Y Y,X\> + \<\na_Y Y, \na_X X\> + 2Sic(\na_Y Y, X) + 2 (\na_YSic)(Y,X) - (\na_X Sic)(Y,Y) \nonumber
\eea 
By the first variational formula, we also have 
\bea \label{naYY,X}
	2\sqrt{\bar \tau}\< \na_Y Y(\bar \tau), X(\bar \tau) \> &=&  \int_0^{\bar\tau} 2 \frac{d}{d\tau} \Big( \sqrt{\tau} \< \na_Y Y, X\> \Big) d\tau \nonumber \\
	 &=& \int_0^{\bar\tau} \sqrt{\tau} \left[ \frac{1}{\tau} \<\na_Y Y, X\> + 2\frac{d}{d\tau}\< \na_Y Y, X \> \right] d\tau  \nonumber \\
	 &=& \int_0^{\bar\tau} \sqrt{\tau} \Big[ 2 \<\na_X \na_Y Y, X\> +(\na_Y Y) S + 4 (\na_Y Sic)(Y,X) - 2(\na_X Sic)(Y,Y) \Big] d\tau \nonumber \\
	 && + \int_0^{\bar\tau} \sqrt{\tau} \<\na_Y Y, 2\na_X X  -\na S + \frac{1}{\tau} X + 4 Sic(X, \cdot) \>d\tau 
\eea 
Thus we obtain
\bea 
Q(Y,Y) &=& \int_0^{\bar\tau} \sqrt{ \tau } \Big[ \na^2S(Y,Y)+2\<R(Y,X)Y, X\> +2 |\na_X Y|^2  \nonumber \\
		&&\quad  -4(\na_Y Sic)(Y,X) + 2 (\na_XSic)(Y,Y) \Big] d\tau
\eea 
As in (\ref{naYY,X}), we have 
\bea \label{naXY,X}
	2\sqrt{\bar \tau}\< \na_X Y(\bar \tau), Y(\bar \tau) \> &=&  \int_0^{\bar\tau} 2 \frac{d}{d\tau} \Big( \sqrt{\tau} \< \na_X Y, Y\> \Big) d\tau \nonumber \\
	 &=& \int_0^{\bar\tau} \sqrt{\tau} \left[ \frac{1}{\tau} \<\na_X Y, Y \> + 2\frac{d}{d\tau}\< \na_X Y, Y \> \right] d\tau  \nonumber \\
	 &=& \int_0^{\bar\tau} \sqrt{\tau} \Big[ \frac{1}{\tau} \<\na_X Y, Y \>+ 4Sic(\na_X Y, Y) +2|\na_X Y|^2  \nonumber \\
	 && + 2 \<\na_X \na_X Y, Y\>  +2(\na_X Sic)(Y,Y) \Big] d\tau
\eea
Hence, $Q(Y,Y)$ can be expressed as
\be
Q(Y,Y) = \int_0^{\bar \tau}2 \< Y, TY \> d\tau + 2 \sqrt{\bar \tau} \< \na_X Y(\bar \tau),Y(\bar \tau)\>.
\ee 
where the vector field $TY$ was defined in (\ref{TY}).

It can be verified that an extension of $Q(Y,Y)$ to a symmetric bilinear form $Q(Y,Z)$ for general vector fields $Y$ and $Z$ is provided by
\be 
Q(Y,Z) = \int_0^{\bar \tau} 2 \sqrt{\tau} \<Z,TY\> + 2\sqrt{\bar \tau}\< \na_X Y(\bar \tau) , Z(\bar \tau)\>.
\ee

\subsection{The reduced distance $L(q,\tau)$ and reduced volume $V(q,\tau)$}

Given $p\in M$, and $t_0 \in (0,T)$, for any $q\in M$ and $0< \bar \tau\leq t_0$, define $L(q,\bar \tau)$ to be the infimal $\cL$-length of all curves $\g$ with $\g(0) =p$ and $\g(\bar \tau) = q$, i.e., 
\be 
	L_{p,t_0}(q,\bar \tau) = \inf _{\g} \cL(\g) =
	 \inf_{\g} \left\{  
		 \int_{0}^{\bar \tau}\sqrt{\tau} \Big( S(\g(\tau)) + |\g'(\tau)|^2 \Big) d \tau : \g(0) = p, \g(\bar \tau) =q 
	 \right\}
\ee
Define the reduced distance based at $(p,t_0)$ by 
\be 
\label{ell}
	\ell_{p,t_0}(q, \tau) = \frac{L_{p,t_0}(q,  \tau)}{2 \sqrt{\tau}},
\ee
and the reduced volume by 
\be
\label{V} 
	\ti V_{t_0}(\tau)  = \int_M (4\pi \tau)^{-\frac{n}{2}} e^{-\ell_{p,t_0}(q,\tau)} dV_{g(\tau)}(q)
\ee 
From now on we will omit the sub-script $p,t_0$ if it is clear from the context. We now compute the derivatives of $L$ and $\ti V$. 

\begin{lemma}
\label{firstderivatives} The first derivatives of $L$ are given by 
\be 
	L_{\bar \tau}(q, \bar \tau) =  2 \sqrt{ \bar \tau }S(q) - \frac{1}{2\bar \tau } L(q, \bar \tau) + \frac{1}{\bar \tau } K
\ee 
and 
\be 
	|\na L|^2(q, \bar \tau) = -4 \bar \tau S(q) +  \frac{2}{\sqrt{\bar \tau}} L(q, \bar \tau) - \frac{4}{\sqrt{\bar \tau}} K
\ee 
where 
\be 
\label{K}
K = \int_0^{\bar \tau}  \tau ^{\frac{3}{2}} H(X(\tau)) d\tau. 
\ee
and $H(X)$ is defined by 
\be 
\label{H}
H(X) = -S_{\tau} -\frac{1}{\tau} S - 2\< \na S, X\> + 2 Sic(X,X)
\ee
\end{lemma}

{\noindent}
{\it Proof.} 
Assume that $q$ is not in the $\bar \tau$ $\cL$-cut locus of $p$ and $\g :[0,\bar \tau]$ be the unique minimizing $\cL$ geodesic from $p$ to $q$ with $\cL$-length $L(q,\bar \tau)$. Let $c:(-\ve, \ve) \to M$ be a curve such that $c(0) = q$, let $\ti \g (u, \tau)$ be the $\cL$ geodesic jointing $p$ and $c(u)$, with $\ti \g (u, \bar \tau) = c(u)$. Let $X(\tau) = \g'(\tau)$ and 
$Y(\tau) = \frac{\p \ti \g}{\p u} \Big|_{u=0}(\tau)$. 
By the first variation of $L$, we have 
\be 
	\<\na L, Y(\bar \tau) \> = \< \na L, c'(0) \> = \frac{d}{du} \Big|_{u=0} L(c(u), \bar \tau)  = 2\sqrt{{\bar \tau}} \< X(\bar \tau) , Y(\bar \tau) \>.
\ee 
hence 
\be
	\na L(q, \bar \tau) = 2 \sqrt{\bar \tau} X(\bar \tau)
\ee 
and 
\be 
	|\na L|^2(q, \bar \tau) = -4 \bar \tau S(q) + 4 \bar \tau  \Big( S(q) + |X(\bar \tau)|^2 \Big). 
\ee 
If we extend the $\cL$-geodesic in $\bar \tau$, we get 
\be 
	\frac{d L(\g(\bar \tau),\bar \tau)}{d \bar \tau} = \sqrt{\bar \tau} \Big( S(\g(\bar \tau)) + |X(\bar \tau)|^2 \Big) 
\ee 
Since 
\be 
	\frac{d L(\g(\bar \tau),\bar \tau)}{d \bar \tau} = L_{\bar \tau}(q, \bar \tau) + \< \na L(q,\bar \tau), X(\bar \tau) \>,
\ee 
we obtain 
\be 
	L_{\bar \tau}(q, \bar \tau) = 2 \sqrt{ \bar \tau }S(q) - \sqrt{\bar \tau} \Big( S(q) +|X(\bar \tau)|^2 \Big).
\ee 

It remains to determine $S(\g(\tau)) +|X(\tau)| ^2$.
First we consider its derivative
\be 
	\frac{d}{d\tau }  \Big( S(\g(\tau)) +|X(\tau)| ^2 \Big) = S_{\tau} + \< \na S, X \> + 2\< \na_X X, X \> + 2 Sic(X,X)
\ee 
Using the geodesic equation \ref{L-geod}, we get 
\be\label{ddtau.S+|X|^2} 
	\frac{d}{d\tau }  \Big( S(\g(\tau)) +|X(\tau)| ^2 \Big) = S_{\tau} + 2 \< \na S, X \>  - \frac{1}{\tau } |X|^2 -2 Sic(X,X) =-H(X) -\frac{1}{\tau} \Big( S(\g(\tau)) +|X(\tau)|^2 \Big) 
\ee  
where $H(X)$ is defined in (\ref{H}).
Multiply both side of (\ref{ddtau.S+|X|^2}) by $\tau^{\frac{3}{2}}$ and integrate from $0$ to $\bar \tau$, we obtain
\be 
	\int_0^{\bar \tau} \tau^{\frac{3}{2}} \frac{d}{d\tau} \Big( S(\g(\tau)) +|X(\tau)| ^2 \Big) d\tau = -K - L(q, \bar \tau)
\ee
where $K$ is defined in (\ref{K}). 
On the other hand, if we integrate the left hand side by parts, we obtain 
\be 
	LHS = \tau^{\frac{3}{2}} \big( S(\g(\tau)) +|X(\tau)|^2 \big) \Big|_0 ^{\bar \tau} - \frac{3}{2} \int_{0}^{\bar \tau} \tau^{\frac{1}{2}}  \Big( S(\g(\tau)) +|X(\tau)| ^2 \Big) d\tau
\ee 
Hence, we have 
\be 
	 \bar \tau ^{\frac{3}{2}} \Big( S(\g(\bar \tau)) + |X(\bar \tau)|^2 \Big) = - K + \frac{1}{2} L(q, \bar \tau)
\ee 
Substituting into the earlier formulas for $L_{\bar\tau}$ and $|\nabla L|^2$ gives the desired formulas. Q.E.D.

\medskip
Next, we consider the second derivatives of $L$. 

\begin{lemma} 
\label{secondderivatives}The following inequality holds in the barrier sense
\be 
L_{\bar \tau}(q,\bar \tau) + \D L(q,\bar \tau) \leq \frac{n}{\sqrt{\bar \tau}}- \frac{1}{2\bar \tau}L(q,\bar \tau),
\ee 
\end{lemma}

\noindent
{\it Proof}. Assume first that $L$ is smooth at $(q,\bar \tau)$.
For a given vector $w\in T_{q} M$, consider the geodesic $c(u)$ such that $c(0) = q, c'(0) =w $. Let $\g_u$ be the unique geodesic from $(p,0)$ to $(c(u),\bar \tau)$,  then the Hessian of the distance function $L$ is given by
\bea
\na^2L(w,w)(q, \bar \tau) = \frac{d^2}{du^2}\Bigg|_{u=0} \cL(\g_u) = Q(Y,Y)
\eea
where $Y$ is the variational vector field of a family of geodesics $\g_u$, hence a $\cL$-Jacobi field with $Y(0) = 0, Y(\bar \tau) = w.$

\smallskip

Our first claim is that for any vector field $W$ with $W(0) = 0, W(\bar \tau) =w$, we have 
$$Q(Y,Y) \leq Q(W,W).$$
Indeed, write $W=Y+Z$, where $Z(0) = Z(\bar \tau) =0$. Therefore, we can find a proper variation $\eta_u$ (of $\g$) with variational vector field $Z$. 
By the minimizing property of $\g$, we have
$$Q(Z,Z) = \frac{d^2}{d u^2}\cL(\eta_u)\Big|_{u=0} \geq 0.$$
On the other hand, by the bilinearity of $Q$, we have
$$Q(W,W) = Q(Y+Z,Y+Z) = Q(Y,Y) + 2Q(Y,Z) +Q(Z,Z)$$
Since $Y$ is a $\cL$-Jacobi field, $TY=0.$ We also have $Z(\bar \tau) =0$, hence $Q(Y,Z) =0$. The claim follows now from the fact that $Q(Z,Z)\geq 0$.

\medskip
Let $w=W(\bar \tau)$ be now a unit vector in $T_{q}M$, and solve for a vector field $W(\tau)$ on $(0,\bar \tau]$ by 
\be 
	\na_{X} W = -Sic(W, \cdot) + \frac{1}{2\tau}W.
\ee 
Then we have 
\be 
	\frac{d}{d\tau}\< W, W\> = 2Sic(W,W) + 2 \< \na_X W, W \> = \frac{1}{\tau} \< W, W \>.
\ee  
hence $|W|^2(\tau) = \frac{\tau}{\bar \tau}.$ Therefore, we can extend $W(\tau)$ continuously by $Y(0) =0$ to $[0,\bar{\tau}]$. We can plug this $W$ into $Q(W,W)$ and get 
\bea 
	Q(W,W) &=& \int_0^{\bar\tau} \sqrt{ \tau } \Big[ \na^2S(W,W)+2\<R(W,X)W, X\> +2 |\na_X W|^2  \nonumber \\
			&&\quad  -4(\na_W Sic)(W,X) + 2 (\na_XSic)(W,W) \Big] d\tau \nonumber \\
			&=& \int_0^{\bar\tau} \sqrt{ \tau } \Big[ \na^2S(W,W)+2\<R(W,X)W, X\> +2 |-Sic(W,\cdot)+ \frac{1}{2\tau} W|^2  \nonumber \\
			&&\quad  -4(\na_Y Sic)(W,X) + 2 (\na_XSic)(W,W) \Big] d\tau \nonumber \\
			&=& \int_0^{\bar \tau} \sqrt{\tau} \Big\{ \na^2{S}(W,W) + 2R(X,W,X,W)  - 4(\na_W Sic)(W,X) + 2(\na_XSic)(W,W) \nonumber \\
			&& + 2 |Sic(W,\cdot)|^2 - \frac{2}{\tau} S(W,W) + \frac{1}{2\tau \bar \tau} \Big\} d\tau 
\eea 
On the other hand, consider $2\sqrt{\bar \tau} Sic(W(\bar \tau), W(\bar \tau))$. 
Note that 
\bea
	\frac{d}{d \tau} Sic(W,W) 
	&=&Sic_{\tau}(W,W) +(\na_X Sic)(W,W) + 2 Sic(\na_X W, W)  \nonumber \\
	&=& Sic_{\tau}(W,W) +(\na_X Sic)(W,W) + \frac{1}{\tau}Sic(W,W) - 2 |Sic(W,\cdot)|^2,
\eea 
and hence
\bea 
	&&2\sqrt{\bar \tau} Sic(W(\bar \tau), W(\bar \tau)) = \int_0^{\bar \tau} 2 \frac{d}{d\tau} \Big( \sqrt{\tau} Sic(W,W) \Big) d \tau  \nonumber \\
	&&\qquad =\int_0^{\bar \tau} \sqrt{\tau} \Big\{ \frac{1}{\tau}Sic(W,W) +Sic_{\bar \tau}(W,W)  +2 (\na_XSic)(W,W)  + \frac{2}{\tau}Sic(W,W) - 4|Sic(W,\cdot)|^2 \Big\} d\tau\nonumber 
\eea 
Adding these two equalities, we get 
\bea 
	&&Q(W,W)+ 2\sqrt{{\bar \tau}} Sic(W(\bar \tau),W(\bar \tau)) - \frac{1}{\sqrt{\bar \tau}} \nonumber \\
	&=& \int_0^{\bar \tau} \sqrt{\tau} \Big\{ \na^2S(W,W) + 2R(X,W,X,W) - 4(\na_W Sic)(W,X) + 4(\na_XSic)(W,W) \nonumber \\
	&& \qquad -2 \big| Sic(W,\cdot)\big|^2 + \frac{1}{\tau} S(W,W) + 2 S_{\tau}(W,W) \Big\} d\tau 
\eea 
We can write this in an easier way 
\be
\label{QWW} 
Q(W,W) = \frac{1}{\sqrt {\bar \tau}}-2\sqrt{\bar \tau}Sic(W(\bar \tau), W(\bar \tau)) - \int_0^{\bar \tau} \sqrt{\bar \tau} H(X,W) d\tau 
\ee
where $-H(X,W)$ is defined to be the integrand above. Note that
if $e_i(\tau)$ is an orthonormal basis of $T_{\g(\tau)}M$, then we have 
$\s_{i=1}^n H(X,e_i) = H(X)$ where $H$ was defined in (\ref{H}).
Since $|W|^2 = \frac{\tau}{\bar \tau}$, we can take $W_i=\Big( \frac{\tau}{\bar \tau}\Big)^{\frac{1}{2}} e_i$, substitute this in the inequality (\ref{QWW}), and obtain
\be 
	\D L \leq \s_{i=1}^n Q(W_i,W_i) = \frac{n}{\sqrt{\bar \tau}} - 2\sqrt{\bar \tau} S - \int_0^{\bar \tau} \sqrt{ \tau} \s_{i=1}^n H(X, W_i) d\tau, 
\ee 
which yields readily
\bea
\D L(q, \bar \tau) \leq 
\frac{n}{\sqrt{\bar \tau}} - 2 \sqrt{\bar{\tau}}S(q,\bar \tau) - \frac{1}{\bar \tau}  K
\eea 
where $K$ was defined in (\ref{K}). Combining this inequality with Lemma \ref{firstderivatives} gives the desired inequality.

\medskip
We now consider the general case.
Let $(q,\bar \tau)$ be any point in the spacetime. If $L$ is not smooth at $(q,\bar \tau)$, then 
$(q,\bar \tau)$ is in the cut-locus of $(p,0)$. Let $\g_1$ be a minimizing geodesic joining $(p,0)$ and $(q,\bar \tau)$ ($\g_1$ does not have to be unique). Given $\ve >0$ small, consider the following barrier function 
$$L_{\ve}(q',\bar \tau') =\int_0^{\ve} \sqrt{\tau}\Big( S(\g_1(\tau)) + |\g_1'(\tau)|^2\Big) d\tau + \inf_{\g} \int_{\ve}^{\bar \tau'}  \sqrt{\tau}\Big( S(\g(\tau)) + |\g'(\tau)|^2\Big) d\tau 
$$
where the infimum is taken over all curves $\g$ with $\g(\ve) = \g_1(\ve)$ and $\g(\bar \tau') = q'.$ It's clear that $L(q',\bar \tau') \leq L_{\ve}(q', \bar \tau')$.  We claim that $L_{\ve}(q',\bar \tau')$ is smooth at $(q,\bar \tau)$. Otherwise, $(q,\bar \tau)$ will be in the cut-locus of $(\g_1(\ve),\ve)$. Then either $(\g_1(\ve),\ve)$ is a conjugate point of $(q,\bar \tau)$ or there exist at least two different minimizing geodesics joining $(\g_1(\ve),\ve)$ and $(q,\bar \tau)$. In both cases, $\g_1$ fails to be minimizing between $(p,0)$ and $(q,\bar \tau).$ 

By definition of the desired inequality in the barrier sense, we need to show that  
$$\Big(\frac{\p}{\p \bar \tau'}+ \D\Big) L_{\ve}(q',\bar \tau') \leq \frac{n}{\sqrt{\bar \tau'}} - \frac{1}{2\bar \tau'}L(q',\bar \tau') + C \ve$$
But we have by the previous calculation
$$\Big(\frac{\p}{\p \bar \tau'}+ \D\Big) L_{\ve}(q',\bar \tau') \leq \frac{n}{\sqrt{\bar \tau'}} - \frac{1}{2\bar \tau} \Big(L_{\ve}(q',\bar \tau') - L(\g_1(\ve),\ve)\Big) \leq \frac{n}{\sqrt{\bar \tau'}} - \frac{1}{2\bar \tau'}L(q',\bar \tau')+C\ve$$
and hence the desired inequality. The proof of Lemma \ref{secondderivatives} is complete.

\medskip
As a simple consequence, the function 
$\bar L(q,\tau) = 2\sqrt{\tau} L(q,\tau)$ satisfies the following inequality in the barrier sense 
\be
\label{Linequality}
	\bar L_{\tau} +\D \bar L \leq 2n,
\ee
Recall that the reduced distance $\ell(q,\tau)$ was defined in (\ref{ell}). It follows immediately
from Lemmas \ref{firstderivatives} and \ref{secondderivatives} that it satisfies the following differential inequalities in the barrier sense,
\bea
	\ell_{\tau} &=& S - \frac{1}{\tau}\ell+ \frac{1}{2\tau^{3/2}}K 
	\label{t derivative eqn}\\
	|\na \ell|^2&=& -S + \frac{1}{\tau}\ell - \frac{1}{\tau^{3/2}}K \label{gradient eqn} \\
	\D \ell &\leq& \frac{n}{2\tau} - S - \frac{1}{2\tau^{3/2}}K \label{Laplacian eqn}
\eea 
which imply in turn
\be
\label{ellinequality} 
	 \ell_{\tau} - \D \ell +  |\na \ell|^2  -S + \frac{n}{2\tau} \geq 0
\ee
Note also that, by the maximum principle applied to (\ref{Linequality}, we have $\min \bar L(\cdot, \tau) \leq 2 n \tau$, and hence 
\be \label{min.l-distance}
\min \ell(\cdot, \tau) = \min \frac{\bar L(\cdot, \tau)}{4\tau} \leq \frac{n}{2}.
\ee 

\subsection{Monotonicity of the reduced volume $\ti V$}

Now we are ready to show that the reduced volume $\ti V(\tau)$ is monotonically nonincreasing in $\tau$. If we let 
$v = (4\pi \tau)^{-n/2} e^{-\ell}$, then we have 
\bea 
	\frac{v_\tau}{v} =- \frac{n}{2\tau} - \ell_{\tau},
	\qquad
	\frac{\D v}{v} = -\D \ell +|\na \ell|^2	
\eea
Hence $v_\tau - \D v + S v \leq 0$. 
Since $\ti V(\tau) = \int_M (4\pi \tau)^{-n/2} e^{-\ell(q,\tau)} dV = \int_M v dV$, we have 
\bea 
	\frac{d}{d\tau} \ti V(\tau)=   \int_M (v_{\tau}+ Sv) dV  \leq \int_M \D v dV  =0
\eea 
i.e., $\ti V(\tau)$ is monotonically nonincreasing in $\tau$.



\section{A $\kappa$-noncollapsing theorem}
\setcounter{equation}{0}

The $\kappa$-noncollapsing theorem of Perelman \cite{P} has been generalized to the coupled Ricci flow by List \cite{L} and M\"uller
\cite{M}. However, we include a statement and proof here for completeness, since the key properties of the reduced distance and volume had already been established in the previous section. First, we recall Perelman's definition of $\kappa$-noncollapse:

\begin{definition}[$\kappa$-noncollapsing]
We say a solution $(M,g,\f)$ of the coupled Ricci flow (\ref{RFHE}) on an interval $[0,T)$ to be $\kappa$-noncollapsing at the scale $\r$, if for each $r<\r$, and all $(x_0,t_0)\in M\ times [0,T)$, the following holds: if $|Rm|(x,t)< r^{-2}$ for all $(x,t)$ in the parabolic neighborhood of $(x_0,t_0)$, $(x,t) \in P(x_0,t_0,r) = B_{t_0}(x_0,r) \times [t_0 -r^2, t_0],$ then $vol(B_{t_0}(x_0,r)) \geq \kappa r^n.$
\end{definition}

\begin{theorem}
Given $n,K,\r,c$ there exists $\kappa = \kappa(n,K,\r,c)$ such that for any solution $(M,g(t),\f(t))$ to Ricci flow coupled with heat equation on $M\times[0,T)$, $T<\infty$  with$(M,g(0))$ complete with $|Rm|+|\na^2 \f|^2 \leq K$ and $inj(M,g(0)) \geq c>0$, then the solution is $\kappa$-noncollapsing at the scale $\r$. 
\end{theorem}

\noindent
{\it Proof.} The proof is by contradiction. Suppose that we can find a sequence of coupled Ricci flows $(M_k,g_k,\f_k)$ such that $|R_m|+|\na^2 \f|^2 \leq K$ on $(M_k,g_k(0))$ and $inj(M_k,g_k(0)) \geq c>0$, and a sequence of $(x_k,t_k)$ and $r_k$ so that $|R_m| \leq r_k^{-2}$ on $B_{t_k}(x_k,r_k) \times [t_k - r_k^2, t_k]$ but 
\be 
	\ve^n_k = \frac{Vol(B_{t_k}(x_k,r_k))}{r_k^n} \to 0.
\ee 

Now, for each $k$, denoting $\tau = t_k - t$, we can define the reduced distance and volume as above. By the short time estimate on curvature and $M$ has bounded geometry at time $t=0$, we know that there exists some $\bar{t}$ such that $M$ has uniformly bounded geometry on $[0,\bar t]$, hence in particular, we can  assume that $t_k \geq \bar t$. 

The idea is to consider $\ti V(\ve_k r_k^2)$ and $\ti V(t_k)$ and show that 
$$0<C^{-1} \leq \lim_{k\to \i}\ti V(t_k) \leq \lim_{k\to \i} \ti V(\ve_k r_k^2) = 0,$$
which is a contradiction. 
The inequality in the middle follows from the non-increasing property of $\ti V$ and the fact that for $k$ large $\ve_k \to 0$. Thus we need to show that we have both
$\lim\limits_{k\to \i} \ti V(\ve_k r_k^2) = 0$, and $\ti V(t_k)$ has a strictly positive lower bound.

\medskip

Note that $\tau = \ve_k r_k^2$ corresponds to time $t = t_k - \ve_k r_k^2$, which is very close to $t_k$ when $k$ large. Let $\g(\tau)$ be a $\cL$-geodesic $\g(\tau)$ with $\g(0) = p$ and initial vector $v= \lim_{\tau \to 0} \sqrt{\tau}X(\tau)$, where $X(\tau) =\g'(\tau)$. 
We claim that if $|v| < \frac{1}{10} \ve_k^{-1/2}$, then the $\cL$-geodesic will stay inside $B_k$ in time $\ve_k r_k^2$. 
Indeed, along the $\cL$-geodesic, we have the geodesic equation (\ref{L-geod}), 
$$\na_X X - \frac{1}{2} \na S + \frac{1}{2\tau} X + 2 Sic(X, \cdot) = 0$$
hence 
\bea
	\frac{d}{d\tau} |X(\tau)|^2 &=& 2Sic(X,X) + 2 \< X, \na_X X \>=
	-\frac{1}{\tau}|X|^2 + \< X,\na S\> -2 Sic(X,X)
\eea 
so 
\be 
	\frac{d}{d\tau}(\tau |X|^2) = -2 \tau Sic(X,X) + \tau\<X, \na S\>
\ee 
By the curvature estimate, we have 
\be 
	|Sic|(x,t)\leq C r_k^{-2} \mbox{ and } |\na S|(x,t)< C r_k^{-3}
\ee 
for $(x,t) \in B(x_k,r_k/2) \times [t_k -r_k^2/2,t_k]$. Hence we have 
\be 
	\Big| \frac{d}{d(\tau/(\ve_k r_k^2))}( \ve_k^{1/2}\sqrt{\tau} |X|) \Big| \leq C\ve_k (\ve_k^{1/2}\sqrt{\tau}|X| ) +C \ve_k^2 \left(\frac{\tau}{\ve_k r_k^2}\right)^{1/2}
\ee  
In particular, on the interval $\frac{\tau}{\ve_k r_k^2} \in [0,1]$ with initial condition 
$\lim_{\tau \to 0} \ve_k^{1/2}\sqrt{\tau} |X| < \frac{1}{10}$,
we have $\ve_k^{1/2}\sqrt{\tau} |X| \leq \frac{1}{9}$ for all $\tau \in [0,\ve_k r_k^2]$ for large $k$, and hence 
\be 
	\int_0^{\ve_k r_k^2} |X(\tau)| d\tau \leq \int_0^{\ve_k r_k^2} \frac{1}{9 \sqrt{\tau}} \ve_k^{-1/2} d\tau = \frac{2}{9}r_k
\ee
From the flow equation that $g_{\tau} = 2Sic$, it follows that the metrics $g(\tau)$ between $\tau=0$ and $\tau = \ve_k r_k^2$ are $e^{C\ve_k}$-biLipschitz close to each other. Thus
\be
	d_{g(t_k)}(x_k,\g(\ve_k r_k^2)) = \int_0^{\ve_k r_k^2} |X(\tau)|_{g(0)} d\tau \leq e^{C\ve} \int_0^{\ve_k r_k}   |X(\tau)|_{g(\tau)} d\tau \leq  \frac{1}{3}r_k
\ee  
for $k$ large enough. Hence the contribution of $\ti V(\ve_k r_k^2)$ coming from those $v\in T_{p_k} M_k$ such that $|v| \leq \frac{1}{10}\ve_k^{-1/2}$ is at most 
$\int_{B_k} (4\pi \ve_k r_k^2)^{-n/2} e^{-\ell(q,\ve_k r_k^2)} dq.$ 

We derive now a lower bound for $\ell(q,\ve_k r_k^2)$ on $B_k$. For $q\in B_k$ assume  $\g:[0,\ve_k r_k^2] \to M$ is the $\cL$-geodesic with $\g(0) = p_k$ and $\g(\ve_k r_k^2) = q$, then 
\be
\cL(\g) \geq \int_0^{\ve_k r_k^2} \sqrt{\tau} S(\g(\tau))d\tau \geq - \int_0^{\ve_k r_k^2}\sqrt{\tau} n(n-1) r_k^{-2} d\tau \geq -C \ve_k^{3/2} r_k
\ee 
hence 
\be 
	\ell(q,\ve_k r_k^2) =  \frac{L}{(\ve_k r_k)^{1/2}} \geq -C \ve_k.
\ee 
Now, the contribution from this part is at most 
\be
	\int_{B_k} (4\pi \ve_k r_k^2)^{-n/2} e^{-\ell(q,\ve_k r_k^2)} dq \leq C \ve_k^{-n/2} r_{k}^{-n} e^{C \ve_k} vol_{t_k -\ve_k r_k^2}(B_k) \leq C \ve_k^{n/2}
\ee 
where we used that $g(t_k-\ve_k r_k^2)$ is close to $g(t_k)$ when $k$ large and $vol(B_k) = \ve_k^n r_k^n$.

To estimate the contribution to $\ti V(\ve_k r_k^2)$ coming from those $v\in T_{p_k}M_k$ such that $|v|\geq \frac{1}{10}\ve_k^{-1/2}$, we can use the monotonicity of the integrand in $\tau$. We claim that as $\tau \to 0$, we have 
$$(4\pi\tau)^{-n/2}e^{-\ell(\cL\exp_{\tau}(v),\tau)} J(v,\tau) \to \pi^{-n/2} e^{-|v|^2}.$$ 

\medskip
The proof of the claim follows  the same  line of the arguments in \cite{CCetc}.
Recall that in terms of the parameter $s=\sqrt\tau$ and $\ti\g(s)=\g(\tau(s))=\g(s^2)$, the reduced length is given by
\bea
L(\ti\g(\bar s),\bar s)=\int_0^{\bar s}\Big({1\over 2}|\ti\g'(s)|^2
+2 s^2S(\ti\g(s))\Big) \,ds.
\eea
Therefore, by L'Hospital's rule, we have 
\be \label{ell limit 1}
	\lim_{\bar s \to 0} \ell(\ti\g(\bar s),\bar s) = \lim_{\bar s\to 0} \frac{L(\ti \g(\bar s),\bar s)}{2 \bar s} = \lim_{\bar s\to 0}\frac{1}{2} \left( \frac{1}{2}|\ti \g'(s)|^2 +2 s^2 S(\ti \g(s)) \right) =  |v|^2
\ee 
Let us study now the limit of the Jacobian $J(v,\tau)$ as $\tau$(or $s$) approaches to $0$. Recall, let $v(u)$ be a curve in $T_p M$, with $v(0) =v$ and $v'(0) = e_i$, where $e_i$ is a orthonormal basis for $T_pM$. Solving the geodesic equation with initial value $\ti \g_u'(0) =2v(u)$ gives a family of geodesics $\ti\g_u(s)$. The variation of this family of geodesic gives a Jacobi field, $J_i(s)$. Therefore 
$$d \cL \exp_{s} : T_pM \to T_{\ti \g(s)} M, \quad  e_i \mapsto J_i(s)$$
Let $h_{ij}(v,s)= \< J_i(s), J_j(s)\>$, then we have 
$$J(v,s)^2 =\det (h_{ij} )(v,s).$$
Note that $\g_u'(0) = 2 v(u)$, and hence
$$2\frac{\p v}{\p u}(0) = \frac{\p^2 \ti \g}{\p u \p s}(0,0) = J_i'(0)$$
i.e., $J_i'(0) = 2 e_i$. Let $E_i(s)$ be the parallel transport of $e_i$ along the geodesic. L'Hospital's rule gives 
\be
	\lim_{s\to 0} \frac{J_i(s) - 2 s E_i(s)}{s} = 0
\ee 
So we have 
$$\lim_{s\to 0} s^{-n}\det(\< J_i, J_j\>)  = \det \Big( \< 2E_i, 2E_j\>\Big) = 2^{2n},$$
hence
\be \label{ell limit 2}
	\lim_{s\to 0}s^{-n/2} J(v,s) \to 2^n
\ee
and by the monotonicity, 
$$(4\pi\tau)^{-n/2}e^{-\ell(\cL\exp_{\tau}(v),\tau)} J(v,\tau) \leq \pi^{-n/2} e^{-|v|^2}.$$
Using the change of variable $q = \cL\exp_{\tau} v$, we have
$$\int_{T_{p_k}M_k\backslash B(0,\frac{1}{10}\ve_k^{-1/2})} (4\pi\tau)^{-n/2}e^{-\ell(\cL\exp_{\tau}(v),\tau)} J(v,\tau) dv \leq \int_{T_{p_k}M_k \backslash B(0,\frac{1}{10}\ve_k^{-1/2})}  \pi^{-n/2}	e^{-|v|^2} \leq e^{-\frac{1}{200\ve_k}}.$$
Combining this estimate with the previous contributions, we obtain
\be
	\lim_{k\to \i} \ti V(\ve_k r_k^2) = 0 .
\ee

Next we estimate $\ti V(t_k)$ from below. Choose a point $q_k$ at time $\frac{\bar t}{2}$ such that $\ell(q_k,t_k - \frac{\bar t}{2}) \leq \frac{n}{2}.$ (This can be done since we have shown in (\ref{min.l-distance}) that $\min \ell(\cdot,\tau)\leq \frac{n}{2}$ for any $\tau$.) Consider a curve $\g_1^{(k)};[0,t_k-\bar t/2] \to M$ with $\g_1^{(k)}(0) = p,$ and $\g_1^{(k)}(t_k - \bar t/2) = q_k$ and another curve $\g_2^{(k)}: [t_k - \bar t/2, t_k]$ with initial point $\g_2^{(k)}(t_k-\bar t/2) = q_k$. Note that $M$ has uniformly bounded geometry on $[0,\bar t/2]$. Thus the distance from $(q_k,t_k-\bar t/2)$ to $(q,t_k)$ for any $q$ in a region around $q_k$ is uniformly bounded, and 
$\ell(\cdot, t_k) \leq C$ on some region of $q_k$. Hence, integrating $e^{-\ell(\cdot, t_k)}$ gives a positive lower bound on $\ti V(t_k)$, where we used $t_k\geq \bar t>0$.
As $\ve_kr_k^2\to 0$ as $k\to \i$, we have $\ve_k r_k^2 < t_k$ for large $k$, and by the (nonincreasing) monotonicity of $\ti V(\tau)$, 
$$\ti V(\ve_k r_k^2) \geq \ti V(t_k).$$
This is a contradiction.

\section{Proof of the Pseudo-Locality Theorem}
\setcounter{equation}{0}

Our proof of the pseudo-locality theorem for the coupled Ricci flow follows the lines of Perelman (\cite{P}, see also \cite{KL}). We begin with the analogues for the coupled Ricci flow of several lemmas of Perelman (Lemma 8.3 in \cite{P}) on the time evolution of the distance function and the fundamental solution of the conjugate heat equation.

\subsection{The time-derivative of the distance function}

\begin{lemma}\label{Box.distance}
Let $(M,g(t), \f(t))$ be a complete Ricci flow coupled with the scalar heat equation.
If  $Ric(g(t_0)) \leq (n-1)K$ in $B_{g(t_0)}(x_0,r_0)$, then for any $x\not\in B_{g(t_0)}(x_0,r_0)$, we have for $d(x,t) = d_{g(t)}(x,x_0)$ 
\be
	\ddt d(x,t) \Big|_{t=t_0} - \D_{g(t_0)} d(x,t_0) \geq -(n-1) \left( \frac{2}{3}K r_0 +r_0^{-1} \right) 
\ee 
\end{lemma}

\noindent
{\it Proof.} Let $\g(s), s\in[0,d(x,t_0)]$ be a normal minimal geodesic with respect to $g(t_0)$ joining $x_0$ and $x$. Then
\be 
	\ddt \Big|_{t=t_0} d(x,t) = -\int_{\g} Sic(\g',\g') ds \geq -\int_{\g}Ric(\g',\g')ds
\ee 
where $Sic = Ric - d\f\otimes d\f.$ Note that $Sic$ is bounded above by $Ric$.

At $x_0$, set $E_1= \g'(0)$, and extend $E_1$ to an orthonormal basis $E_i$ of $T_{x_0}M$. Parallel transporting this basis along $\g$ gives us an orthonormal basis of $T_{\g(s)}M$. 

Recall the second variation formula of the distance function. Assume $w\in T_{x}M$, and let $W$ be the Jacobi field given by $W(x_0) = 0$, $W(x) =w$. Then 
\be 
	\na^2{d}(w,w) = I(W,W) \leq I(Y,Y)
\ee 
where $Y$ is any vector field with $Y(x_0) = 0$, $Y(x) = w$ and $I(W,W)$ is the index form defined by
\be 
I(V,W) = \int_{\g}\Big( \<V',W'\>-\<V',\g'\>\<W',\g'\>- R(\g',V,\g',W) \Big) ds.
\ee 
Thus
\be 
	\D_{g(t_0)} d(x,t_0) = \sum_{i=1}^n \na^2{d}(e_i,e_i) \leq \sum_{i=1}^n \int_{\g} |\na_{\g'} F_i|^2 - \left<\na_{\g'}F_i,\g'\right>^2  -R(\g',F_i,\g',F_i) ds
\ee 
for any vector field $F_i$ with $F_i(x_0) =0$ and $F_i(x) = E_i(x)$.
Choosing $F_i$ along  $\gamma$ as
\bea
F_i(\gamma(s)) = 
	\left\{
		\begin{array}{ll}
		\frac{s}{r_0}E_i(s) \quad &s\in[0,r_0] \\
		E_i(s) \quad & s\in [r_0,d(x,t_0)].
		\end{array}
	\right.
\eea 
we obtain 
\be
\< \na_{\g'} F_i, \g'\> = 0 
\ee 
for all $i = 2,\cdots, n$ and $s\in [0,d(x,t_0)]$  and 
$ |\na_{\g'} F_1|^2 - \left<\na_{\g'}F_1,\g'\right>^2 = 0$. Hence 
\bea 
\D_{t_0} d(x,t_0) &\leq & \int_{\g} \sum_{i=2}^n |\na_{\g'} F_i|^2  - R(\g',F_i,\g',F_i) ds \nonumber \\
&=&  \int_0^{r_0} \sum_{i=2}^ n \Big( \frac{1}{r_0^2} - \frac{s^2}{r_0^2}R(\g',E_i,\g',E_i) \Big) ds + \int_{r_0}^{d(x,t_0)}\sum_{i=2}^n -R(\g',E_i,\g',E_i) ds \nonumber  \\
&=& \frac{n-1}{r_0} + \int_0^{r_0} (1-\frac{s^2}{r_0^2}) Ric(\g',\g') ds - \int_0^{d(x,t_0)} Ric(\g',\g') ds \nonumber \\
&\leq& \frac{n-1}{r_0} +(n-1)K(r_0 - \frac{1}{3}r_0) + \ddt d(x,t) \Big|_{t=t_0}\label{Perelman estimate}
\eea  
This completes the proof.

By similar argument, we have the following lemma, for which we omit the proof.
\begin{lemma}
\label{Perelman1}
Let $(M,g(t),\f(t))$ be a complete Ricci flow coupled with the scalar heat equation.
Assume  $Ric(g(t_0),x) \leq (n-1)K$ for any $x\in B_{g(t_0)}(x_1,r_0) \cup B_{g(t_0)}(x_2,r_0)$, and $d_{t_0}(x_1,x_2)\geq 2r_0$.
Then we have 
\be 
	\ddt d_t(x_1,x_2) \Big|_{t=t_0} \geq -2(n-1) \left( \frac{2}{3}Kr_0 + r_0^{-1} \right)
\ee
where $d_t(x_1,x_2) = d_{g(t)}(x_1,x_2).$ 
\end{lemma}

\subsection{The localized $\cW$-functional and conjugate heat equation}

Let $\Box = \p_t - \D$. The formal adjoint $\Box^*$
of $\Box$ is defined by the relation
\bea
 0 = \int_a^b \int_M \Box \varphi \cdot \psi dV dt - \int_a^b \int_M  \varphi\cdot \Box^*\psi dV dt 
\eea
for any $a,b \in [0,T]$, $\varphi, \psi \in C_0^\i(M\times(a,b))$. It is readily recognized that $\Box^*= - \p_t - \D + S$. 

\smallskip

Let $ u = (4\pi(T-t))^{-\frac{n}{2}} e^{-f}$ be a solution of the conjugate heat equation $\Box^* u =0$. Set
\bea
\label{v}
v = \Big( (T-t)(2\D f - |\na f|^2 +S ) + f - n \Big) u .
\eea
Then by direct calculations, we have
\bea
\label{Box}
\Box^* v = -2(T-t) \left(\left| S_{ij} +\na_i \na_j f - \frac{g_{ij}}{2(T-t)} \right|^2 + (\D \f - \< \na \f,\na f\>)^2 \right) u
\eea
and  in particular we have $\Box^* v\le 0$. We remark that the integral of $v$ over $M$ is the $\mathcal W$-functional introduced by Perelman (see e.g. \cite{P,CZ,Ni})
$$\cW(g_{ij}, f, T-t) = (4\pi(T-t))^{-n/2}\int_M\Big( (T-t)(|\nabla f|^2 +S) + f - n  \Big)e^{-f} dV,$$ in case the integration by parts holds.

Similar to the Ricci flow case (see \cite{P,Ni}), we have the following differential Harnack estimate for the fundamental solutions to the conjugate heat equation.
\begin{lemma}
\label{Perelman2}
Let $u$ be the fundamental solution based at $(p,T)$, namely $\Box^* u = 0$ and
$u(x,t) \to \d_p(x)$ as $t\to T$. Then we have 
$v\leq 0$ for all $t<T$, where $v$ is defined as in (\ref{v}). 
\end{lemma}

\noindent
{\it Proof.} For any fixed $t_0\in (0,T)$ and any positive function $h_0$, solve the equation $\Box h = 0$ with initial $h(t_0) = h_0$. Then we have 
$$\frac{d}{dt} \int_M h v dV = \int_M \Big( (\Box h ) v - h (\Box^* v) \Big) dV = -\int_M h(\Box^* v) dV \geq 0.$$
Hence, $\int_M hv dV$ is increasing in $t\in (t_0,T)$. Moreover, as stated by Perelman (\cite{P}, see also \cite{Ni})
\be 
\lim_{t\to T^-} \int_M h(x,t) v(x,t) dV_{g(t)} = 0,
\ee
hence we have 
\be 
\int_Mh(x,t_0) v(x,t_0) dV_{g(t_0)} \leq 0.
\ee
Since $h(x,t_0)$ and $t_0$ are arbitrary, we have $v(x,t) \leq 0$ for any $t<T.$

\subsection{A point selection lemma}

We return now to proof of the pseudo-locality theorem. 
Without loss of generality, we can assume $r_0=1$. 
Assume that the theorem is not true, that is, there exists 
$\ve_k,\d_k \to 0$ such that for each $k$, there exists a complete solution $(M_k,g_k,\f_k,p_k)$ to coupled Ricci flow (\ref{RFHE}), such that $S(g_k(0)) \geq -1$ on ball $B_{g_k(0)}(p_k,1)$,  $Area_{g_k(0)}(\p \O)^n \geq (1-\d_k)c_n Vol_{g_k(0)}(\O)^{n-1}$ and $|\phi_{k,0}|\le C$, for any $\O \subset B_{g_k(0)}(p_k,1)$, and $(x_k,t_k)$ such that $d_{g(t_k)}(x_k,p_k)\leq \ve_k$ but
$$|Rm|(x_k,t_k) \geq \a t_k^{-1}+ \ve_k^{-2} .$$ 
Moreover, we can choose a smaller $\ve_k$ such that 
\be  \label{Rm.bound}
	|Rm|(x,t) \leq \a t^{-1} + 2 \ve_k^{-2}
\ee 
for all $t\in (0,\ve_k^2)$ and $d_{g_k(t)}(x,p_k) \leq \ve_k$. 
We divide the argument into several steps.

\medskip

The first step is to choose some other point $(\bar x,\bar t)$, such that the Riemannian curvature tensor can be controlled by $|Rm|(\bar x,\bar t)$ in a parabolic neighborhood of $(\bar x, \bar t)$, provided there exists an $(x,t)$ satisfying the above hypotheses. 

The following point-selection lemmas can be proved in the same way as the Ricci flow case, so we omit the proof and refer to that of {claim 1} and {claim 2} in section 10 of \cite{P} (see also \cite{KL}).
\begin{lemma}
\label{Perelman} For any large $A>0$ and any solution   $(M,g(t),\f(t),p)$ to coupled Ricci flow, if there is a point $(x_0,t_0) \in M\times (0,\ve^2]$ such that $|Rm|(x_0,t_0) \geq \alpha t_0^{-1} + \ve^{-2}$ and $d_{t_0}(x_0,p)\leq \ve$, then there is a point $(\bar x,\bar t)\in M_{\a}$ such that $d_{\bar t}(\bar x, p) \leq (1+2A)\ve $ and 
\be 
	|Rm|(x,t) \leq 4 |Rm|(\bar x,\bar t);
\ee  
for any $(x,t)\in M_\a$, $0<t\leq \bar t$ such that 
$$
	d_t(x,p) \leq d_{\bar t}(\bar x, p) +A |Rm|^{-1/2}(\bar x,\bar t),
$$
where 
$$M_\a:= \{ (x,t)\in M\times (0,\ve^2] : |Rm|(x,t) > \a t^{-1} \}.$$
\end{lemma}

\begin{lemma}
\label{Perelman4}
Under the same assumption as previous lemma, the point $(\bar x,\bar t)$ selected above satisfies 
\be
	|Rm|(x,t) \leq 4Q =: 4|Rm|(\bar x,\bar t)
\ee 
for any $x\in B_{g(\bar t)}(\bar x, \frac{1}{10} A Q^{-1/2} ) \times [\bar t - \frac{1}{2}\alpha Q^{-1}, \bar t].$ 
\end{lemma}

Applying  Lemma \ref{Perelman4}, we can find $(\bar x_k,\bar t_k) \in M_k \times (0,\ve_k^2]$ with $d_{g_k(\bar t_k)}(\bar x_k, p_k) \leq (1+2A_k)\ve_k$,  satisfying the above properties, i.e.,
\be
	|Rm|_{g_k(t)}(x,t) \leq 4 Q_k = 4|Rm|(\bar x_k,\bar t_k)
\ee 
for any $(x,t)\in \O_k$, where $\Omega_k$ is defined by
\bea
\label{Omega}
\Omega_k=\{(x,t);\ d_{\bar t_k}(x,\bar x_k) \leq \frac{1}{10}A_k Q_k^{-1},
\ \bar t_k- \frac{1}{2}\alpha Q_k^{-1} \leq t \leq \bar t_k\}
\eea 

\subsection{A gap lemma}

For each $k$, let $u_k$ be the fundamental solution at $(\bar x_k,\bar t_k)$
of the conjugate heat equation,
i.e., $u_k$ is the solution of
$$\Box_t^* u_k = 0$$
with initial condition Dirac function $\d_{\bar x_k}(x)$ at time $\bar t_k$. Define as in (\ref{v})
$$v_k = \Big((\bar t_k - t)(S+ 2 \D f_k - |\na f_k|^2) +f_k -n \Big)u_k $$
where we have set $u_k = \Big( 4\pi(\bar t_k - t) \Big)^{-n/2} e^{-f_k}.$  

\begin{lemma}
\label{gap}
There exists a uniform constant $b>0$(independent of $k$) and a $\ti t_k \in  (t_k- \frac{1}{2}\alpha Q_k^{-1}, \bar t_k)$ for each $k$, such that 
\be
	\int_{B_k} v_k dV_{g_k(\ti t_k)} \leq -b < 0
\ee 
where $B_k = B_{g_k(\ti t_k)}(\bar x_k, \sqrt{\bar t_k - \ti t_k}).$ 
\end{lemma}

\noindent
{\it Proof.} The proof is by contradiction. Assume that it is not true, which means 
that for any $\ti t_k \in  (t_k- \frac{1}{2}\alpha Q_k^{-1},  \bar t_k)$, there exists a subsequence, still denoted by $k$, such that
\be \label{int.vk.nonneg} 
 \liminf_{k \to \i} \int_{B_k} v_k dV_{g_k(\ti t_k)}  \geq 0 .
\ee 
Consider the following rescaling, 
\be
\hat{g}_k(t) = Q_k g_k(Q_k^{-1} t + \bar t_k), \hat{\f}_k(t) = \f_k(Q_k^{-1}t +\bar t_k)
\ee
for $t\in [-Q_k \bar t_k, 0]$. Then $(\hat g_k(t), \hat\f_k(t))$ also satisfies the coupled Ricci flow (\ref{RFHE}). Note that under the parabolic rescaling
$\O_k$ becomes
\be 
	\O_k = B_{\hat{g}(0)}(\bar x_k, \frac{1}{10}A_k) \times [ -\half\a ,0] 
\ee 

Now, we consider the following two cases: either along a subsequence the injective radius of $\hat g_k(0)$ at $\bar x_k$ has positive lower bound, or there is no such a lower bound along any subsequence. 

\medskip

{\bf Case 1}: by List's compactness theorem (Theorem 7.5 in \cite{L}), we can find a subsequence, again denoted by $k$, such that $(M_k,\hat{g}_k(t), \hat\f_k(t), p_k)$ converges in the pointed $C^\i$-CG (Cheeger-Gromov) sense to a new complete coupled Ricci flow $(M_\i, g_{\i}, \f_{\i},p_{\i})$ for $t\in [ -\half \a, 0]$, and $|Rm|_{g_\i}(x,t) \leq 4$ for all $(x,t)\in M_\i \times [-\half\a,0]$ and $|Rm|(x_\i,0) = 1$. 
For each $M_k$, we have the fundamental solution $\hat u_k$ based at $(\bar x_k,\bar t_k)$, $\hat u_k$ converge to $u_\i$ in the same sense, where $u_\i$ is a fundamental solution to the conjugate heat equation on $(M_\i, g_\i, \f_\i, x_\i)$ based at $(x_\i,0)$. So $v_k$ converge to $v_\i$ in the same sense, and $v_\i \leq 0$ by Lemma \ref{Perelman2}. By the assumption (\ref{int.vk.nonneg}), for any fixed $t_0\in [-\half \a, 0]$, we have 
\be 
\int_{B_{g_\i(t_0)}(x_\i, \sqrt{-t_0})} v_\i(\cdot ,t_0) dV_{g_\i(t_0)} \geq 0,
\ee
Thus $v_\i(\cdot ,t_0) = 0$ on $B_{g_\i(t_0)}(x_\i, \sqrt{-t_0})$. Next we
show that
\be 
	v_\i \equiv 0 \quad \mbox{on } M_\i \times (t_0,0]
\ee 
For this, take a non-negative, non-trivial, function $h_0$ with compact support in $B_{g_\i(t_0)}(x_\i, \sqrt{-t_0})$, and solve the following heat equation 
\be 
	\p_t h = \D_{g_\infty(t)} h, \quad h(\cdot,t_0) = h_0(\cdot)
\ee 
By the maximum principle, we have $h(x,t)>0$ for any $t>t_0$ and 
\be 
	\frac{d}{dt} \int_{M_\i} h(x,t) v_\i(x,t) dV_{g_\i(t)} = \int_{M_\i} \Box_t h  v_\i - h \Box_t^* v_\i \geq 0.
\ee 
Since the integral $\int h v dV$ is zero at both $t=t_0$ and $t=0$ (see \cite{P,Ni}), it is zero for any $t\in (t_0,0)$. Since $h$ is strictly positive, and $v_\i$ non-positive, our assertion follows.

\medskip 

The formula (\ref{Box}) for $\Box_t^* v_\i$ shows that $(g_\i(t),\f_\i(t))$ is a complete extended Ricci soliton, i.e.,
\bea\label{Ricci soliton}
Sic_\i(g_\i(t))+ \na^2 f_\i(t) - \frac{g_\i(t)}{-2t} &=& 0 \\
\D \f_\i - \< \na \f_\i, \na f_\i \> &=& 0\label{Ricci soliton 2}
\eea 
for all $t\in (t_0,0)$. We require the following version for the coupled Ricci flow of a result of Zhang \cite{Z}, 

\begin{theorem}
\label{coupled Zhang}
Let $(M_\infty,g_\infty(t),\phi_\infty,f_\infty)$ be as in (\ref{Ricci soliton}) and (\ref{Ricci soliton 2}),
then $S_\infty(t)\ge 0$ and  the gradient vector $-\nabla_{g_\infty} f_\infty $ is a complete vector field, i.e. it generates a one-parameter family of diffeomorphisms $\varphi(t): M_\infty \to M_\infty$ for all $t\in (-\infty, 0)$.
\end{theorem}

The proof of Theorem \ref{coupled Zhang} will be given in the Appendix.

We can now apply Theorem \ref{coupled Zhang}.
The completeness of $g_\i(t)$ implies that of $\na f_\i$. Thus the vector field $-\na f_\i$ can be integrated to give a family of diffeomorphisms $\varphi(t): M_\i \to M_\i$, $t_0\le t<0$, such that 
\bea 
\ddt \varphi(t) =- \na f_\i(t)\circ \varphi(t)  \nonumber
\eea
Consider $\ti g(t) = -t^{-1}\varphi(t)^* g_\i(t)$. We have 
\bea
\frac{\partial}{\partial t}\ti g(t) &=& t^{-2}\varphi(t)^* g_\i(t) - t^{-1} \varphi(t)^* (L_{-\na f_\i}g_\i(t) - 2 Sic(g_\i(t) )  \nonumber\\
&=& - 2 t^{-1} \varphi(t)^* \left( -Sic(g_\i(t)) -\na^2 f_\i(t) - \frac{g_\i(t)}{2t}  \right) =0.
\eea 
Hence $\ti g(t)$ is independent of $t\in (t_0, 0)$, the Riemannian curvature of $\ti g(t)$ is bounded and $|Rm_{\tilde g}|(x_\i) \not = 0$.  On the other hand, 
$|Rm|_{\ti g(t)} =(- t) |Rm|_{g_\i(t)}$, therefore $|Rm|_{g_\i(t)}(x_\i) = \frac{|Rm|_{\ti g(t)}(x_\i)}{-t} \gg 1$ when $t$ is close to $0$. This contradicts $|Rm|_{g_\i(0)}(x_\i) = 1$.

\medskip
{\bf Case 2}:
Suppose now that there is a subsequence so that the injectivity radii of the metrics at $\bar x_k$ tend to zero, w.l.o.g., assume 
$r_k = inj(\bar x_k, \hat{g}_k(0)) \to 0$. Rescale $(M_k,\hat{g}_k, \hat{\f}_k,p_k)$ by 
\be 
\ti g_k(t) = r_k^{-2} \hat{g}_k(r_k^2t), \quad \ti \f_k(t) = \hat \f_k(r_k^2 t).
\ee 
where $t\in [-\half \a r_k^{-2} , 0]$. The region $\O_k$ becomes
\be
	d_{\ti g_k(0)}(x,\bar x_k) \leq \frac{1}{10}A_k r_k^{-1} \to \infty, \quad  t\in  [-\half \a r_k^{-2} , 0]
\ee 
The injectivity radius of $\ti g_k(0)$ at $\bar x_k$ is $1$. On $\O_k$, the Riemannian curvature tensor is bounded above by $4r_k^2$, hence we get a subsequence converging in $C^\i$-CG sense to a complete coupled Ricci flow $(M_\i, g_\i, \f_\i, p_\i)$ for $t\in (-\i, 0]$. Moreover the uniform curvature bound on $\O_k$ implies that the solution $(M_\i, g_\i)$ is a flat metric. Hence, by similar arguments as in the first case, we get a family of solitons $(M_\infty,g_\infty,\phi_\infty,f_\infty)$ satisfying (\ref{Ricci soliton}) and (\ref{Ricci soliton 2}). Since $Rm(g_\infty(t)) = 0$, by Theorem \ref{coupled Zhang}, $S(g_\infty) = - |\nabla_{g_\infty}\phi_\infty|^2\ge 0$, we have $\phi_\infty = $const. Therefore
\be 
	\na^2 f_\i  = \frac{g_\i(t)}{-2t}.
\ee 
By the uniformization theorem in Riemannian geometry, the universal cover of $(M_\i,g_\i)$ is isometric to $(\R^n,g_{can})$, $\pi:{ \mathbf {R}}^n\to M_\infty$. Pulling back to the universal cover, we get 
$$\na^2 \pi^*f_\i = \frac{g_{can}}{-2t}> 0$$
Hence, $\pi^* f_\i$ is a strictly convex function in $\R^n$. Since a convex function on $\R^n$ can never be periodic, $\pi$ has to be trivial. Therefore, we have $(M_\i, g_\i(t)) = (\R^n, c(t)g_{can})$, but this contradicts the fact that $inj_{g_\i}(x_\i,0)$ is finite. The proof of Lemma \ref{gap} is complete.

\bigskip

We establish next another version of the preceding gap lemma, but with the volume form $dV_{g_k(\ti t_k)}$ replaced by $dV_{g_k(0)}$. For this, we need the assumption on the initial metrics, which we had not used as yet. 

\smallskip

Let $\varphi$ be a smooth function on $\R$ which is one on $(-\i,1]$, decreases to zero on $[1,2]$, and is zero on $[1,\infty)$, with 
$(\varphi')^2 \leq 10\varphi$ and $-\varphi'' \leq 10\varphi.$

The following construction is done for each individual $(M_k,g_k,\f_k,p_k)$. Let $\ti d_k(x,t) = d_k(x,t) + 200 n \sqrt{t}$ where $d_k(x,t)=d_{g_k(t)}(p_k,x)$, and define a function $h_k(x,t)$ by 
$$h_k(x,t) = \varphi\left(\frac{\ti d_k(x,t)}{10 A_k\ve_k}\right)$$

\begin{lemma}
\label{gap1}
If the constants $A_k$ are chosen to be large enough,
then we have for all $k$
\bea
\int_M h_kv_k dV_{g_k}\Big|_{t=0}\leq -{1\over 2}b<0.
\eea
where $b$ is the constant in Lemma \ref{gap}.
\end{lemma}

\medskip
\noindent
{\it Proof.}
We suppress the subindex $k$ for notational simplicity. It is easy to show that 
\be 
	\Box_t h(x,t) = \frac{1}{10 A \ve }\left( d_t -\D d + \frac{100n}{\sqrt{t}} \right) \varphi'\left(\frac{\ti d(x,t)}{10 A \ve }\right)  - \frac{1}{(10A\ve)^2}\varphi''\left(\frac{\ti d(x,t)}{10 A \ve }\right)
\ee 
We claim that, 
if $A$ is large enough, then 
\bea
\label{Boxh}
\Box_t h(x,t) \leq  - \frac{1}{(10A\ve)^2}\varphi''\left(\frac{\ti d(x,t)}{10 A \ve }\right).
\eea
For this, it suffices to show that if $A$ is large enough, then 
on the support of $\varphi'(\ti d(x,t)(10 A \ve)^{-1})$, i.e., when $10 A\ve \leq \ti d(x,t) \leq 20 A \ve $, we have 
\be 
	d_t - \D d + \frac{100n}{\sqrt{t}} \geq 0.
\ee
Indeed, for $t\in [0,\ve^2]$, we have $200n\sqrt{t} \leq A\ve$ if $A$ is large enough, and hence 
$$9A\ve \leq d(x,t) \leq 21 A\ve $$
Let $r_0 = \sqrt{t}.$
Since $r_0 \leq \ve$, we obviously have $x \not \in B(p,r_0)$. Moreover, from (\ref{Rm.bound}), we know that 
$|Rm|(x,t) \leq \a t^{-1} + 2 \ve^{-2}$ for $x \in B(p,r_0)$. Thus we can apply Lemma \ref{Box.distance} to get 
\bea 
d_t - \D d &\geq& - (n-1) \left( \frac{2}{3} (\a t^{-1} + 2\ve^{-2}) t^{1/2}+ t^{-1/2} \right)  \nonumber \\
&=& -(n-1) \left(\frac{2}{3}\a + \frac{4}{3}\ve^{-2} t + 1 \right) t^{-1/2}  
\geq -100 n t^{-1/2}.
\eea 
This establishes the claim.


\medskip
Recall that $\Box_t ^* v \leq 0$. Thus 
\bea
\frac{d}{dt} \left( \int_M h (-v) dV \right) &=& \int_M\Big( \Box h (-v) + h \Box^* v \Big) dV  \leq  \int_M \Box h (-v) dV \nonumber \\
&\leq & \frac{1}{(10A\ve)^2} \int_M  - \varphi''(-v) dV \nonumber \leq \frac{1}{(10A\ve)^2} \int_M 10\varphi (-v) dV \\
&\leq & \frac{10}{10 (A\ve)^2} \int_M h (-v) dV 
\eea 
and hence 
\be 
 \frac{d}{dt} \log \int_M h(-v) dV \leq \frac{1}{10 (A\ve)^2}.
\ee 
Integrating from $0$ to $\ti t$ yields  (here $\tilde t$ is from Lemma \ref{gap})
\be
\frac{\int_M h (-v) dV\Big|_{t= \ti t}} {\int_M h (-v) dV\Big|_{t= 0} }  \leq e^{\frac{\ti t}{10 A^2\ve^2}} \leq e^{\frac{1}{10 A^2}} 
\ee
where we used $\ti t \leq \ve^2.$ 
Hence, we get 
\be
\int_M h(-v) dV \Big|_{t=0}  \geq e^{-\frac{1}{10A^2}} \int_M h(-v) dV |_{t = \ti t}.
\ee 
It suffices to show that $h\equiv 1$ on the ball $B = 
B_{g(\ti t)}(\bar x, \sqrt{\bar t - \ti t})$, and then the desired inequality follows from Lemma \ref{gap}. 

To see this, it suffices to show for any $x\in B_{g(\tilde t)}(\bar x, \sqrt{\bar t - \tilde t})$, it holds that
$d_{g(\tilde t)}(x,p)\le 3 A\varepsilon$. Recall by Lemma \ref{Perelman4},
$|Rm|(y,t)\le 4 Q = 4 |Rm|(\bar x,\bar t),$ for all $$(y,t)\in B_{g(\bar t)}\big(\bar x, \frac {1}{10} A Q^{-1/2}\big)\times [\bar t - \frac 1 2 \alpha Q^{-1}, \bar t].$$ And by distance comparison argument we have
$$d_{g(\tilde t)}(y,\bar x)\le e^{4Q(\bar t - \tilde t)} d_{g(\bar t)}(y,\bar x)\le e^{2\alpha} d_{g(\bar t)}(y,\bar x) .$$
For $x\in B_{g(\tilde t)}(\bar x, \sqrt{\bar t - \tilde t})$, we have
\bea d_{g(\tilde t)}(x,p)&\le d_{g(\tilde t)}(x,\bar x) + d_{g(\tilde t)}(p,\bar x)\nonumber\\
& \le \sqrt{\bar t - \tilde t} + e^{2\alpha} d_{g(\bar t)}(p,\bar x)\nonumber\\
&\le \sqrt{\frac{\alpha Q^{-1}}{2}} + e^{2\alpha} (1+2A)\varepsilon\nonumber\\
& \le \varepsilon + (1+2A)\epsilon\le 3A\varepsilon,\nonumber
\eea
if $A$ is large enough and we use the estimate  $Q\ge \alpha \bar t^{-1}\ge \alpha \epsilon^{-2}$. Hence the desired statement follows.

\bigskip
We show next that the gap inequality in Lemma \ref{gap1} can be expressed in a form closer to that appears in log-Sobolev inequalities.

We continue to suppress the subindex $k$ for notational simplicity.
Set $u = (4\pi(\bar t - t))^{-n/2} e^{-f}$, and $v = \Big( (\bar t -t)(S + 2\D f -|\na f|^2) + f - n \Big) u$. Then we can write 
\bea
\int_M h(-v) dV\Big|_{t=0}
= 
\int_M -  \Big( \bar t(S + 2\D f -|\na f|^2) + f - n  \Big) (4\pi \bar t)^{-n/2} he^{-f}  dV 
\eea
In terms of $\ti u \equiv u h$, $\ti f \equiv f - \log h$, 
this can be rewritten as 
\bea
\int_M h(-v) dV\Big|_{t=0}
&=&
\int_M - \Big( \bar t (S + 2 \D \ti f + 2 \D \log h - |\na\ti f + \na \log h|^2) + \ti f + \log h - n \Big) (4\pi \bar t)^{-n/2} e^{-\ti f} dV \nonumber \\
&=& \int_M \Big( \bar t(- 2\D \ti f +|\na \ti f|^2) - \ti f +n\Big) (4\pi \bar t)^{-n/2} e^{\ti f} dV \nonumber \\
&& + \int_M \Big(-\bar t S + \bar t |\na \log h|^2 -  \log h \Big) (4\pi \bar t)^{-n/2} e^{-\ti f} dV
\eea
after an integration by parts. 

The second term in the above right hand side can be estimated as follows.
At time $t = 0$, we have $S(g,\f)(0) = R(g) - |\na \f|^2 \geq -1$. Thus
\be 
	\int_M -\bar t S \ti u dV \leq \int_M \bar t \ti u dV 
	\leq \int_M \bar t u dV\leq \bar t \leq \ve^2.
\ee
From the definition of $h$, it follows immediately that 
$$|\na \log h |^2 = \frac{1}{(10 A\ve)^2} \frac{\varphi'^2}{\varphi^2} \leq \frac{1}{(10A\ve)^2} \frac{10}{\varphi}
$$
Since $1\leq\varphi\leq 2$ on the support of $\varphi'$, we have then
\be
	\int_M \bar t |\na \log h|^2 h u dV \leq \frac{\bar t}{10 A\ve^2} \int_M u dV \leq \frac{1}{10 A^2}
\ee
For the last term $-\int_M h \log h u dV$, by the definition of $h$, it's non-zero only on the region $B(p, 20 A\ve) \backslash B(p, 10 A\ve).$ Since $ - x \log x \leq \frac{1}{e} \leq 1 $ on $(0,1]$, we obtain
\be 
\int_M -(h\log h)\, u dV \leq \int_{B(p,20 A\ve) \backslash B(p,10 A\ve)} u dV \leq \int_M u dV - \int_{B(x_0,10 A\ve)} u dV
\ee 
Finally, a similar argument as above shows that
$$\int_{B(p,10A\ve)} u dV \geq \int_M \varphi\Big({d(x,t)\over 5A\e}\Big)\, u dV \geq 1 - cA^{-2},$$
for some constant $c$. Putting all these together, and applying Lemma \ref{gap1},
we obtain
\be 
\label{gap3}
\int_M (-\bar t |\na \ti f|^2 - \ti f + n) (4\pi \bar t)^{-n/2} e^{-\ti f}dV \geq (1-A^{-2}) {b} - (1+c) A^{-2} -\ve^2\ge \frac b 2.
\ee 
if $A$ is large and $\varepsilon$ is small enough.

Define $\ti g = \frac{1}{2\bar t} g$, then $d\ti V = (2\bar t)^{-n/2} dV$, $\bar t|\na \ti f|_{g}^2 =  \frac{1}{2}|\na \ti f|^2_{\ti g}$. Now we restore the subscript $k$, and normalize the function $ \ti u_k $ as follows
\be 
\label{Uk}
	U_k = \frac{\ti u_k}{\int_{M_k} \ti u_k d \ti V_{\ti g_k}} = \frac{(2\bar t_k)^{n/2} \ti u_k}{\int_{M_k} \ti u dV_{g_k} }.
\ee
Define then the function $F_k$ by 
\be
\label{Fk}
	U_k = (2\pi)^{-n/2} e^{-F_k}.
\ee 

\begin{lemma}
\label{gap4}
Let $F_k$ and $U_k$ be defined as above. Then we have for $k$ large
\be \label{Contradict.log.Sobolev}
	\int_M \Big( -\half |\na F_k|^2 - F_k + n \Big) U_k d\ti V_{\ti g_k} \geq \frac{1}{2}b > 0.
\ee 
\end{lemma}

\noindent
{\it Proof}. Write
$\frac{\ti u_k}{\si_k} = (4\pi \bar t_k)^{-n/2} e^{-F_k}$, where $\si_k \equiv \int_{M_k} \ti u_k dV_{g_k} \to 1, \mbox{ as } k \to \i$. From the definition of $\ti f_k$, we have $F_k = \ti f_k + \log \si_k$. Plugging this to the inequality (\ref{gap3}), we get
\be 
	\int_{M_k} \Big( -\frac{1}{2} |\na F_k|^2 - F_k + \log \si_k + n \Big) U_k \si_k d\ti V_{\ti g_k} \geq (1- A_k^{-2}) b -(1+c) A_k^{-2} - \ve_k^2
\ee  
Since $A_k \to \i$, $\ve_k \to 0$ and $\si_k \to 1$ as $k\to \i$, the lemma follows.

\subsection{Logarithmic Sobolev inequalities}

We now derive the desired contradiction. Once the estimate in Lemma \ref{gap4} is established, this part of the argument is the same as for the Ricci flow. Nevertheless, we give it for the sake of completeness.

\medskip

We begin by recalling the log Sobolev inequality on $\R^n$. One version of it says that
\be 
	\int_{\R^n} \Big( - \frac{1}{2}|\na F|^2 - F +n \Big) (2\pi)^{-n/2} e^{-F} dV \leq 0
\ee 
for compactly supported functions $U = (2\pi)^{-n/2} e^{-F}$ satisfying $\int_{\R^n} U dV = 1$. It would have been easy to compare this inequality with the inequality established in Lemma \ref{gap4}, had both been on $\R^n$, but we need a version of the log Sobolev inequality on Riemannian manifolds. For this, it is convenient to make use of another version for $\R^n$, which can be obtained as follows.

\smallskip

Let $U_c(x) = c^n U(cx)$, and define $F_c(x)$ by $U_c(x) = (2\pi)^{-n/2} e^{-F_c(x)}$. Then $F_c(x) = F(cx) - n \log c$. Applying the logarithmic Sobolev inequality to $U_c(x)$ ( since $\int_{\R^n} U_c(x)dx = 1$) gives
\be 
	\int_{\R^n} \Big(-\frac{1}{2}c^2 |\na F(cx)|^2 - F(cx) + n\log c + n \Big) (2\pi)^{-n/2}e^{-F(cx)} c^n dx \leq 0
\ee   
Change of variable by $y = cx$, and maximize the right hand side term for $c \in (0, \i)$. A simple calculation shows that the maximum is achieved at 
$$c^2 = n(\int_{\R^n}|\na F(x)|^2 (2\pi)^{-n/2} e^{-F(x)} dx )^{-1}$$
Substituting $c$ back into the inequality above, we obtain the inequality
\be 
	\int_{\R^n} |\na F|^2U dx \geq n \exp \Big( 1 - \frac{2}{n} \int_{\R^n} FU dx \Big) 
\ee 
for all $U = (2\pi)^{-n/2} e^{-F}$ and $\int_{\R^n} Udx = 1$.

\medskip

Using symmetrization arguments, we show now that the same inequality holds on a Riemannian manifold, provided the isoperimetric inequality holds for domains on that Riemannian manifold.  More precisely,

\begin{lemma}
Let $(M,g_{ij})$ be a compact Riemannian manifold. Assume that 
$Area(\p \O)^{n} \geq (1-\d) c_n Vol(\O)^{n-1}$ for any $\O \subset B_g(p,r)$, then for any $U$ with compact support on $B_g(p,r)$, with $\int_{M} U dV_g= 1$, we have 
\be \label{Riem.log.Sobolev}
	\int_{M} |\na F |^2 U dV_g \geq (1- \d)^{2/n} n \exp \Big(1 - \frac{2}{n} \int_{M} F U dV_g\Big)
\ee 
\end{lemma}

\noindent
{\it Proof.} We use another equivalent formulation of the log Sobolev inequality. On $\R^n$, it says that
\be \label{Euclid.log.Sobolev.form2}
	\int_{\R^n} 4 |\na \varphi|^2 dx \geq n \exp \left( 1 + \log 2\pi + \frac{2}{n}\int 2(\log \varphi)\, \varphi^2 dx \right)
\ee 
and the Riemannian version, equivalent to (\ref{Riem.log.Sobolev}), is the following
\be \label{Riem.log.Sobolev.form2}
	\int_{M_k} 4|\na \varphi|^2 dV_g \geq (1-\d)^{2/n} n \exp \left( 1 + \log2\pi + \frac{2}{n}\int_{M_k} 2(\log{\varphi}) \varphi^2  \right)
\ee 
where we have set $\varphi^2 = U = (2\pi)^{-n/2} e^{-F}$, $\int_{M} \varphi^2 dV_g= 1$ and $F = -2\log \varphi - \frac{n}{2}\log 2\pi.$ 

\medskip

To prove (\ref{Riem.log.Sobolev.form2}), define $\phi^*$ on $(\R^n,g_{can})$ as the symmetrization of $\varphi$ as follows. More specifically, let $M(s) = \{  x\in M: \varphi(x) \geq s\} $, and definee $\phi^*$ to be the non-increasing, radially symmetric function on $\R^n$ satisfying $\phi^*(0)=\sup_{M} \varphi$, and for any $s$,   
$$Vol_{g_{can}}(\{ \phi^*(x) \geq s\}) = Vol_{g}(M(s)).$$
Recall the co-area formula states that for any function $g$, we have 
\be 
	\int_{\O} g(x) |\na u| dx = \int_{-\i}^{\i}\int_{u^{-1}(s)} g(x) d\si_s ds
\ee 
Applying this to $M(s)$, with $g(x) = \frac{1}{|\na \varphi|}$, and differentiating with respect to $s$, we obtain 
$$Vol(M(s))' = \int_{\varphi^{-1}(s)} \frac{1}{|\na \varphi|} d\si= \int_{(\phi^*)^{-1}(s)} \frac{1}{|\na \phi^*(s)|} d\si_{g_{can}}$$
On the other hand, since $\phi^*$ is symmetric, we have $|\na \phi^*|$ is constant of $\{ \phi^* = s\}$. Hence, 
\bea 
Area( \p\{ \phi^* = s\} )^2 &=& \int_{\{ \phi^* = s\}} |\na \phi^*| d\si_{g_{can}} \int_{\{ \phi^* = s\}} \frac{1}{|\na \phi^*|} d\si_{g_{can}} \nonumber \\
&=& \big(Area( \p\{ \phi^* = s\} )^n \big)^{2/n}  = \big( c_n Vol(\{\phi^*\geq s\})^{n-1} \Big)^{2/n} \nonumber \\ 
&=& \big( c_n Vol(M(s))^{n-1} \big)^{2/n} \nonumber \\
&\leq& (1-\d)^{-2/n} Area^2(\p M(s)) \nonumber \\
&\leq& (1-\d)^{-2/n} \int_{\{ \varphi = s\}} |\na \varphi| d\si_g \int_{\{ \varphi = s\}} \frac{1}{|\na \varphi|} d\si_g.
\eea 
Hence, we get 
\be 
\int_{\{ \varphi = s\}} |\na \varphi| d\si_g\geq (1-\d)^{2/n} \int_{\{ \phi^* = s\}} |\na \phi^*| d\si_{g_{can}}
\ee 
Applying the co-area formula again, we get 
\be 
	\int_{M} |\na \varphi|^2 d V_g  \geq (1-\d)^{2/n} \int_{\R^n} |\na \phi^* |^2 d x,
\ee
and 
\bea 
	\int_{M} \varphi^2 \log \varphi  dV_g&=& \int_{0}^\i \int_{\{\varphi = s\}} \frac{\varphi^2 \log \varphi}{|\na \varphi|} d\si_g\nonumber ds\\
	&=& \int_0^\i s^2 \log s \int_{\{\varphi= s\}} \frac{1}{|\na \varphi|}d\si_gds \nonumber \\
	&=& \int_0^\i s^2 \log s \int_{\{\phi^*= s\}} \frac{1}{|\na \phi^*|}d\si_{g_{can}} ds \nonumber \\
	&=& \int_{\R^n} \phi^*{}^2 \log \phi^* dx
\eea 
Putting all these together gives the logarithmic inequality (\ref{Riem.log.Sobolev.form2}) on Riemannian manifolds.

\subsection{Completion of the proof of the pseudo-locality theorem}

We return now to the setting of $\ti g = \frac{1}{2\bar t} g$, then $d\ti V = (2\bar t)^{-n/2} dV$, with the functions $U_k$ and $F_k$ defined as in (\ref{Uk}) and (\ref{Fk}). 
The equations (\ref{Contradict.log.Sobolev}) and (\ref{Riem.log.Sobolev}) imply that 
\be 
 (1-\d_k)^{2/n} \frac{n}{2} \exp \left( 1 - \frac{2}{n}\int_{M_k} F_k U_k d\ti V_k \right) + \int_{M_k} F_k U_k d\ti V_k 	- n  \leq -\frac{b}{2}
\ee
Let $\eta_k = 1 - \frac{2}{n}\int_{M_k} F_k U_k d\ti V_k$, the left hand side term can be written as 
$$
	LHS = \frac{n}{2} \left( (1-\d_k)^{2/n} e^{\eta_k} - 1 -\eta_k \right)	\leq -\frac{b}{2}
$$
However, consider the function $(1-\d_k)^{2/n} e^x - 1 -x $. Its minimum occurs at $x_0 $ given by $ (1-\d_k)^{2/n} e^{x_0}  = 1$, and hence its minimum is $ - x_0 =  \frac{2}{n} \log(1-\d_k)$. Thus, we have 
\be 
	LHS \geq \frac{n}{2}\cdot  \frac{2}{n}\log(1-\d_k) = \log(1-\d_k) \to 0, \mbox{ as } k \to \i.
\ee 
This is a contradiction, and the proof of the pseudo-locality theorem is complete.

\section{Convergence of parabolic rescaling to a soliton}
\setcounter{equation}{0}

The goal of this section is to establish Theorem \ref{2}. It is shown in the first subsections that,  in presence of a Type I singularity, the parabolic rescalings of the coupled Ricci flow will converge to a soliton. The non-triviality of the soliton will be established in the last subsection by applying the pseudo-locality theorem.

\subsection{Estimates for the reduced distance}

Suppose the flow satisfies the Type I condition (\ref{Type I}), for any $T'<T$, we have $\sup_M |Rm|(g(t))\le \frac{C_0}{T'-t}$ for any $t\in [0,T')$. For any fixed $t_0\in [0,T')$, take $r^2 = T'-t_0$, and by the Type I assumption
$$\sup_{t\in [0,t_0], x\in B_{g(t)}(p,r)} r^2|Rm|(x,g(t))\le \sup_{t\in [0,t_0]}(T'-t_0) \frac{C_0}{T' - t} = C_0,$$
thus by the derivative estimates (\ref{List derivative estimate}), it follows that
$$\sup_{B_{g(t)}(p,r/2)}|\nabla^m Rm|(g(t))\le C(n,m)C_0^{1+\frac m 2} \bk{ \frac{1}{T'-t_0} + \frac{1}{t} }^{1+\frac m 2},\quad \forall t\in (0,t_0].$$
In particular, at $(p,t_0)$, we have
\begin{equation}\label{higher order estimate}
|\nabla Rm|(p,t_0)\le \frac{C(n,C_0)}{(T'-t_0)^{3/2}},\quad |\nabla^2 Rm|(p,t_0)\le \frac{C(n,C_0)}{(T'-t_0)^2},\quad |\partial_t Rm|(p,t_0)\le \frac{C(n,C_0)}{(T'-t_0)^2},
\end{equation}
and since $p\in M,t_0\in [0,T')$ are arbitrary, the above estimates hold on $M\times [0,T')$.

\medskip

Set $\tau=\tau(t) = T'-t$. We will denote $\bar g(\tau) = g(T'-\tau)$. By the higher order estimates (\ref{higher order estimate}), we have
\begin{equation}\label{estimates for K}
\max(|S|_{\bar g(\tau)}, |Sic|_{\bar g(\tau)})\le \frac{C(n,C_0)}{\tau},\quad \max( \abs{\nabla S}_{\bar g(\tau)}, |\partial_\tau S|_{\bar g(\tau)}  )\le \frac{C(n,C_0)}{\tau^2}.
\end{equation}
We are now ready to establish the following estimates for the reduced distance, which are obtained by Naber (\cite{N}) in the Ricci flow case.

\begin{lemma}\label{lemma estimate}
		There exists a constant $C(n,C_0)$ depending only on the Type I constant $C_0$ such that
		\begin{enumerate}
			\item $C(n,C_0)^{-1}\frac{d^2_{\bar g(\bar \tau)}(p,q)}{\bar \tau} - C(n,C_0) \le \ell(q,\bar \tau)\le C(n,C_0) \frac{d^2_{\bar g(\bar \tau)}(p,q)}{\bar \tau} + C(n,C_0)$.
			\item $\abs{\nabla \ell}\le \frac{C(n,C_0)}{\bar \tau} \bk{ 1+ \frac{d^2_{\bar g(\bar\tau)}(p,q)}{\bar \tau} }$.
			\item $\Big| \frac{\partial}{\partial \tau} \ell   \Big|\le \frac{C(n,C_0)}{\bar\tau} \bk{ \frac{d^2_{\bar g(\bar \tau)}(p,q)}{\bar \tau}  + 1 }$,
		\end{enumerate}
where $\ell$ is the reduced distance with base space-time $(p,T')$.
	\end{lemma}

\noindent
{\it Proof.} 	We apply the estimates in (\ref{gradient eqn}). For this, we estimate the quantity $K$ in (\ref{K}) as follows,
\bea
	|K|&\le& C(n,C_0) \int_0^{\bar \tau}\tau^{3/2}\bk{\frac{1}{\tau^2 } + \frac{|X|}{\tau^{3/2}} + \frac{\abs{X}}{\tau}  } d\tau \nonumber \\
	&\le& C(n,C_0)\int_0^{\bar \tau} \tau^{3/2}\bk{\frac{1}{\tau^2} + \frac{\abs{X}}{\tau}    } d\tau \nonumber \\
	&\le& C(n,C_0)\int_0^{\bar \tau} \sqrt{\tau}( S + \abs{X}  ) d\tau + C(n,C_0) \sqrt{\bar \tau}\nonumber \\
	& = &C(n,C_0) \sqrt{\bar \tau}\ell + C(n,C_0) \sqrt{\bar \tau}, \nonumber
\eea
Substituting the above estimates into (\ref{gradient eqn}), we get the estimate
\begin{equation}\label{gradient 1}
\abs{\nabla \ell}\le \frac{C(n,C_0)}{\bar \tau} \ell + C(n,C_0).
\end{equation}




Next, we claim that there exists a constant $C_2 = C_2(n,C_0)$ such that 
	$$|\ell(p,\bar \tau)|\le C_2,\quad \forall \bar \tau \in (0,T'].$$
Indeed, consider the constant map $\gamma:[0,\bar \tau]\to M$ defined by $\gamma(\tau)\equiv p$. The $\mathcal L$-length of this curve is bounded above by 
	$$\cL(\gamma)\le \int_0^{\bar \tau}\sqrt{\tau} \frac{C(n,C_0)}{\tau}d\tau \le C(n,C_0)\sqrt{\bar \tau},$$
 then the claim follows at once from the definition of $\ell(p,\bar \tau)$.

\medskip
Thus from (\ref{gradient 1}) we have
\begin{equation}\label{gradient 2}
|\nabla \sqrt{\ell + C_1}|_{\bar g(\bar \tau)}\le \frac{C(n,C_0)}{\sqrt{\bar \tau}},
\end{equation}
The mean value theorem implies
$$\sqrt{\ell(q,\bar \tau) +C_1} - \sqrt{\ell(p,\bar \tau) + C_1}\le C(n,C_0) \frac{d_{\bar g(\bar \tau)}(p,q)}{\sqrt{\bar \tau}},$$
which is 
\begin{equation}\label{gradient 3}
\ell(q,\bar \tau)\le C(n,C_0) \frac{d^2_{\bar g(\bar \tau)}(p,q) }{\bar \tau} + C(n,C_0).
\end{equation}
This is the desired upper bound for $\ell(q,\bar \tau)$.
To derive the  lower bound, we begin by showing that,
for any minimal geodesic $\sigma$ with respect to the metric $\bar g(\bar \tau)$, we have
	$$\int_\sigma Sic(\sigma'(s),\sigma'(s))ds\le \frac{C(n,C_0)}{\sqrt{\bar\tau}}.$$

	If $L_{\bar g(\bar \tau)}(\sigma)\le 2\sqrt{\bar \tau}$, this integral inequality follows directly from the curvature assumption $\sic\le \ric\le \frac{C_0}{\bar \tau}$. 
	
	If $d\equiv L_{\bar g(\bar \tau)}(\sigma)>2\sqrt{\bar \tau}$, take an orthonormal basis $\{E_i\}_{i=1}^n$ of vector fields which are parallel with respect to $\bar g(\bar \tau)$ along the geodesic $\sigma$, and $E_n = \sigma'(s)$.  Choose a function 
	$\phi(s)$ satisfying $\phi(s)={s\over \sqrt {\bar\tau}}$ for $s\in [0,\sqrt{\bar\tau}]$, 
	$\phi(s)=1$ for $s\in[\sqrt{\bar\tau},d-\sqrt{\bar\tau}]$, and $\phi(s)={d-s\over\sqrt{\bar\tau}}$ for $s\in [d-\sqrt{\bar\tau},d]$.
By the second variational formula of the length function (applied to each variational vector field $\phi E_i$ for $1\le i\le n-1$), we have
	$$0\le \int_0^d \sum_{i=1}^{n-1} \abs{\nabla_{E_n} (\phi E_i) } - \int_0^d \phi^2(s) Ric(E_n, E_n) ds.$$
Consequently,
	$$\int_\sigma Ric(\sigma',\sigma')\le \frac{2(n-1)}{\sqrt{\bar \tau}} + \frac{n-1}{\bar \tau} C_0 \int_0^d (1-\phi^2)ds = \frac{2(n-1)}{\sqrt{\bar \tau}} + \frac{4(n-1)C_0}{3\sqrt{\bar \tau}}.$$
	The claim follows now from $Sic\le Ric$.
	
	\medskip

Choose now $\mathcal L$-geodesics $\gamma_1$ and $\gamma_2:[0,\bar \tau]\to M$ connecting the space-time points $(p,0)$ with $(p,\bar \tau)$ and $(q,\bar \tau)$, respectively. Let $\sigma_\tau:[0,d_{\bar g(\tau)}(\gamma_1(\tau), \gamma_2(\tau))]\to M$ be the minimal geodesic connecting the points $\gamma_1(\tau)$ and $\gamma_2(\tau)$ under the metric $\bar g(\tau)$. Then
\bea
	\frac{d}{d\tau} d_{\bar g(\tau)}(\gamma_1(\tau),\gamma_2(\tau)) &=& \innpro{\nabla d, \dot\gamma_1}_{\bar g(\tau)} + \innpro{\nabla d, \dot\gamma_2}_{\bar g(\tau)} + \int_\sigma Sic(\sigma'(s),\sigma'(s))ds \nonumber \\
	&\le&  |\dot\gamma_1|_{\bar g(\tau)} + |\dot\gamma_2(\tau)|_{\bar g(\tau)} + \frac{C(n,C_0)}{\sqrt{\tau}}. \nonumber 
\eea
	Integrating the above inequality from $\tau\in [0,\bar\tau]$ yields
	$$d_{\bar g(\bar \tau)}(p,q)\le \int_0^{\bar \tau} \bk{ |\dot\gamma_1| + |\dot\gamma_2|} d\tau + C(n,C_0)\sqrt{\bar \tau}.$$
The first term on the right hand side can be estimated by
	$$\int_0^{\bar \tau} |\dot\gamma_1|\le \int_0^{\bar \tau}\sqrt{\tau} \abs{\dot\gamma_1} + C(n) \sqrt{\bar \tau}\le 2\sqrt{\bar \tau} \ell(p,\bar \tau) + C(n,C_0) \sqrt{\bar\tau}\le C(n,C_0) \sqrt{\bar \tau},$$
while the second term can be estimated by
	$$\int_0^{\bar \tau} |\dot\gamma_2|d\tau\le \bk{\int_0^{\bar \tau} \sqrt{\tau} \abs{\dot \gamma_1}  }^{1/2} \bk{\int_0^{\bar \tau} \frac{1}{\sqrt{\tau}}  }^{1/2}\le C(n,C_0) \sqrt{\bar \tau}\bk{ \ell(q,\bar\tau) + C(n,C_0)  }^{1/2}, $$
Combining with these estimates we get 
	$$d_{\bar g(\bar \tau)}(p,q)\le \sqrt{\bar \tau}\Big(C(n,C_0) + C(n,C_0)\bk{ \ell(q,\bar \tau) + C(n,C_0)  }^{1/2}\Big),$$
	that is 
	$$C(n,C_0)^{-1} \frac{d^2_{\bar g(\bar \tau)}(p,q)}{\bar \tau} - C(n,C_0) \le \ell(q,\bar \tau).$$
This is the desired lower bound for $\ell(q,\bar\tau)$. The remaining estimates in Lemma \ref{lemma estimate} follow from the equations (\ref{t derivative eqn}) and (\ref{gradient eqn}), so the proof is completed.

\bigskip

We restate the lemma above in terms of the original extended Ricci flow. For any $T'\in (0,T)$, there exists a constant $C=C(n,C_0)$ such that (denote $\ell_{(p,T')}(q,t) = \ell(q, T'-t)$ as in the previous notation)
\begin{equation}\label{reduced 1}
C(n,C_0)^{-1} \frac{d^2_{g(t)}(p,q)}{T'-t} - C(n,C_0)\le \ell_{(p,T')}(q,t) \le C(n,C_0) \frac{d^2_{g(t)}(p,q)}{T'-t} +C(n,C_0),
\end{equation}
\begin{equation}\label{reduced 2}
\abs{\nabla_{g(t)}\ell_{(p,T')}(q,t)}\le \frac{C(n,C_0)}{T'-t}\bk{ \frac{ d^2_{g(t)}(p,q)  }{T'-t} + 1  },
\end{equation}
\begin{equation}\label{reduced 3}
\Big|\frac{\partial}{\partial t} \ell_{(p,T')}(q,t) \Big|\le  \frac{C(n,C_0)}{T'-t}\bk{ \frac{ d^2_{g(t)}(p,q)  }{T'-t} + 1  },
\end{equation}
and the following inequality hold in the distributional sense,
\begin{equation}\label{reduced 4}
\frac{\partial\ell}{\partial t} + \Delta_{g(t)}\ell - \abs{\nabla \ell} + S - \frac{n}{2(T'-t)}\le 0.
\end{equation}

Take a sequence of times $\{T_i\}$ such that $T_i\nearrow T$, and we have the corresponding reduced distance functions $$\ell_i(q,t) = \ell_{(p,T_i)}(q,t): M\times [0,T_i]\to \R,$$
and each $\ell_i$ satisfies the estimates (\ref{reduced 1}), (\ref{reduced 2}), (\ref{reduced 3}) and (\ref{reduced 4}) with the same constant $C(n,C_0)$. In particular on any compact subset $\Omega\subset M\times [0,T)$, the functions $\ell_i$ satisfy uniform $C^0$, $C^1_M$, and $C^1_t$ estimates. In addition, $$\| \ell_i\|_{W^{1,2}(\Omega)}\le C(n,C_0,\Omega).$$
Thus up to a subsequence and a diagonal argument we may assume that
$$\ell_i \rightharpoonup \ell_\infty \in W^{1,2}_{loc}(M\times [0,T))$$
weakly in $W^{1,2}_{loc}(M\times [0,T))$, and on any compact subset $\Omega\subset M\times [0,T)$, the convergence is uniform in $C^0(\Omega)$ norm. Moreover, it is not hard to see the limit function is locally Lipschitz, so the limit $\ell_\infty \in W^{1,2}_{loc}$, and $\ell_\infty$ satisfies similar estimate as $\ell_i$ in (\ref{reduced 1}). 

We define
\bea
V_\infty(t)= (4\pi(T-t))^{-n/2} \int_M e^{-\ell_\infty(q,t)}dV_{g(t)}.
\eea

\begin{lemma}
The function $V_\infty(t)$ satisfies the following properties:

{\rm (1)} $V_\infty(t)\leq 1$, $\forall t\in [0,T)$.

{\rm (2)} $V_\infty(t_1)\leq V_\infty(t_2)$, for $t_1<t_2\in (0,T)$. 

\end{lemma}

\noindent
{\it Proof.} In view of Lemma \ref{lemma estimate}, we have for each $\ell_i$ and $t\in [0,T)$,
$$e^{-\ell_i(q,t)}\le C(n,C_0) e^{- \frac{ d^2_{g(t)}(p,q)  }{C(n,C_0)(T_i - t)}}\le C(n,C_0) e^{- \frac{ d^2_{g(t)}(p,q)  }{C(n,C_0)(T- t)}}$$ 
Thus the function on the right hand side in integrable with respect to $dV_{g(t)}$,
so by the Lebesgue dominated convergence theorem 
\begin{equation}\label{reduced volume convergence}\lim_{i\to \infty} V_i(t) = \lim_{i\to \infty}(4\pi(T_i - t))^{-n/2}\int_M e^{-\ell_i(q,t)}dV_{g(t)}
=V_\infty(t)  \end{equation}
It follows from (\ref{ell limit 1}) and (\ref{ell limit 2}) that $\lim_{t\to T_i} V_i(t) = 1$ (see also \cite{P}) and $V_i(t)$ is nondecreasing in $t$, so for each $t\in (0,T)$ when $i$ is large enough $V_i(t)\le \lim_{s\to T_i}V_i(s) = 1$. Part (1) of the lemma follows.
Similarly, Part (2) follows from the monotonicity of $V_i(t)$ for each $i$, so that $V_i(t_1)\leq V_i(t_2)$ for $t_1<t_2\in (0,T)$. Q.E.D.

\medskip
\begin{lemma}
\label{ell infinity}
The function $\ell_\infty$ satisfies the following inequality  in the distribution sense,
\bea
\p_t\ell_\infty+\Delta_{g(t)}\ell_\infty-|\na\ell_\infty|^2+S-{n\over (T-t)}\geq 0.
\eea
\end{lemma}

\noindent
{\it Proof.} 
Since $\ell_i$ converge weakly to $\ell_\infty$ in $W^{1,2}_{loc}(M\times [0,T))$, for any vector field $V$ on $M$ with compact support and $\psi\in C^\infty_c((0,T))$ we have
\begin{equation}\label{weak convergence}
\int_0^T \int_M \bk{ \psi \innpro{V, \nabla \ell_i}_{g(t)} + \frac{d\psi}{dt} \ell_i  }dV_{g(t)}dt \to \int_0^T \int_M \bk{ \psi \innpro{V, \nabla \ell_\infty}_{g(t)} + \frac{d\psi}{dt} \ell_\infty}dV_{g(t)}dt.
\end{equation}
We aim to prove that for any nonnegative $\varphi\in C^\infty_c(M\times (0,T))$, it holds that
\begin{equation}\label{estimate to go}
\int_0^T \int_M \bk{ -\frac{\partial\varphi}{\partial t} \ell_\infty + \Delta_{g(t)}\varphi \ell_\infty - \varphi \abs{\nabla \ell_\infty} + S \varphi - \frac{n\varphi}{2(T-t)}    }dV_{g(t)}dt\le 0.
\end{equation}
First note that by (\ref{reduced 4}), the following holds for $\ell_i$ for $i$ large enough, 
\begin{equation}\label{estimate to go 1}
\int_0^T \int_M \bk{ -\frac{\partial\varphi}{\partial t} \ell_i + \Delta_{g(t)}\varphi \ell_i - \varphi \abs{\nabla \ell_i} + S \varphi - \frac{n\varphi}{2(T_i-t)}    }dV_{g(t)}dt\le 0.
\end{equation}
By locally $C^0$ uniform convergence of $\ell_i$ to $\ell_\infty$, it follows that the first, second, fourth and last terms in (\ref{estimate to go 1}) converge to those in (\ref{estimate to go}) as $i\to\infty$, respectively. So to show (\ref{estimate to go}), it suffices to show the third terms also converge, i.e.
\begin{equation}\label{to show}
\int_0^T \int_M \varphi \abs{\nabla \ell_i} dV_{g(t)}dt \to \int_0^T \int_M \varphi \abs{\nabla \ell_\infty} dV_{g(t)}dt.
\end{equation}
Since $\ell_\infty$ is the weak limit of $\ell_i$ in $W^{1,2}_{loc}(M\times (0,T))$, it follows that
$$\int_0^T \int_M \varphi \abs{\nabla\ell_\infty} dV_{g(t)}dt \le \liminf_{i\to \infty}\int_M \varphi \abs{\nabla \ell_i}dV_{g(t)}dt.$$
So to show (\ref{to show}) it suffices to prove that 
\begin{equation}\label{eqn sup}\limsup_{i\to \infty}\int_0^T \int_M \varphi \abs{\nabla \ell_i}dV_{g(t)}dt \le \int_0^T \int_M \varphi \abs{\nabla \ell_\infty} dV_{g(t)}dt,\end{equation}
noting that
\bea \label{sup estimate}
\limsup_{i\to \infty}\int_0^T \int_M \varphi \abs{\nabla \ell_i}dV_{g(t)}dt  & =& \limsup_{i\to \infty}\int_0^T \int_M \varphi\innpro{ \nabla (\ell_i - \ell_\infty -\epsilon_i),\nabla \ell_i  }dV_{g(t)}dt\nonumber\\
&&\quad + \limsup_{i\to\infty}\int_M  \innpro{ \varphi \nabla \ell_\infty, \nabla\ell_i  } dV_{g(t)}dt, 
\eea 
using (\ref{weak convergence}), we see the second term on RHS of (\ref{sup estimate}) converges to $\int_0^T \int_M \varphi\abs{\nabla \ell_\infty} dV_{g(t)}dt$. Since $\ell_i$ converges uniformly to $\ell_\infty$ on supp$\varphi\subset M\times (0,T)$, there exists a sequence of numbers $\epsilon_i>0$ which tend to $0$ such that
$$\ell_\infty -\ell_i + \epsilon_i\ge 0,\qquad \mbox{on supp}\,\varphi.$$
To deal with the first term on RHS of (\ref{sup estimate}), by (\ref{Laplacian eqn}) it follows that 
$$\Delta_{g(t)}\ell_i\le \frac{C(n,C_0)}{T_i - t} + C(n,C_0)\frac{\ell_i}{T_i - t},$$
in the distribution sense. Multiplying both sides by $\varphi(\ell_\infty -\ell_i + \epsilon_i)$ and integrating it follows that
$$\int_0^T \int_M \innpro{\nabla(\varphi(\ell_i - \ell_\infty -\epsilon_i)), \nabla \ell_i}\le \int_0^T \int_M \frac{C(n,C_0)}{T_i - t} \varphi(\ell_\infty - \ell_i + \epsilon_i) (C(n,C_0) + \ell_i), $$
the RHS tends to $0$ as $i\to \infty$, because both $\ell_\infty$ and $\ell_i$ are bounded on supp$\,\varphi$ and $\ell_\infty - \ell_i + \epsilon_i\to 0$. The LHS is equal to
$$\int_0^T \int_M \varphi\innpro{\nabla (\ell_i - \ell_\infty -\epsilon_i),\nabla \ell_i} + \int_0^T \int_M (\ell_i - \ell_\infty -\epsilon_i)\innpro{\nabla \varphi,\nabla \ell_i}, $$ and the second term tends to $0$, so $$\limsup_{i\to \infty }\int_0^T \int_M \varphi\innpro{\nabla (\ell_i - \ell_\infty -\epsilon_i),\nabla \ell_i}\le 0.$$ Thus we finish the proof of (\ref{to show}), and also that of  (\ref{estimate to go}). 

\bigskip

Define function $u_\infty$ by
\bea
u_\infty=(4\pi(T-t))^{-n/2}e^{-\ell_\infty}.
\eea
Then the inequality in Lemma \ref{ell infinity} is equivalent to the inequality
\bea
\label{Boxu}
\Box^* u_\i = (-\partial_t - \Delta_{g(t)} + S ) u_\i\le 0
\eea
in the distribution sense.
So it follows that
\bea
\label{monotone V}
\frac{d}{dt} V_\infty = \frac{d}{dt}\int_M u_\infty dV_{g(t)} = \int_M (\partial_t + \Delta _{g(t)} - S ) u_\infty \ge 0.
\eea

\begin{lemma}\label{lemma 15}
If $V_\infty(t_1)=V_\infty(t_2)$ for some $t_1<t_2\in (0,T)$, then $g(t)$ is a coupled gradient soliton, that is,
$$Sic + {\nabla^2}\, \ell_\infty - \frac{g}{2(T-t)} = 0,\quad \Delta \phi - \innpro{\nabla \ell_\infty, \nabla \phi} = 0.$$
\end{lemma}

\noindent
{\it Proof.} The existence of two such values $t_1$ and $t_2$ implies that the integrand on the right hand side of (\ref{monotone V}), which is known to be $\geq 0$ by Lemma \ref{ell infinity}, must vanish identically. By parabolic regularity, 
the function $u_\infty$ is actually $C^\infty$ in $M\times [0,T)$. Thus the function $v_\infty$ defined by
$$v_\infty = \bk{(T-t)(S + 2\Delta \ell_\infty - \abs{\nabla \ell_\infty}) + \ell_\infty - n  }u_\infty 
$$
as well as $\Box^*v_\infty$ must vanish identically. Since we have, by a direct calculation,
$$
\Box^* v_\infty =\bk{-2(T-t)\Big|{Sic + {\nabla^2}\,\ell_\infty - \frac{g}{2(T-t)}}\Big|^2 - 2 \abs{\Delta \phi - \innpro{\nabla \phi,\nabla\ell_\infty}} } u_\infty.
$$
The vanishing of $\Box^*v_\infty$ implies immediately that $g(t)$ is an extended soliton, as claimed.

\subsection{Blow-ups of Type I $\kappa$-noncollapsed solutions}

Let $(M,g(t),\phi(t))$ be an coupled Ricci flow solution to (\ref{RFHE}) with $|\phi_0|\leq C$,  and assume that it is of Type I, i.e., $$\sup_M |Rm|(g(t))\le \frac{C_0}{T-t}, \quad t\in [0,T),$$
and $T<\infty$ is the maximal existence time of the flow. 

Take any sequence of numbers $\lambda_i\to \infty$, and define a sequence of coupled Ricci flows
$$g_i(t) = \lambda_i g(\lambda_i^{-1}t + T),\quad \f_i(t) = \f(\lambda_i^{-1}t + T),\quad \forall t\in [-\lambda_i T, 0).$$
By the Type I condition we have
$$\sup_M |Rm|(g_i(t)) = \lambda_i^{-1}\sup_M |Rm|(g(\lambda_i^{-1}t + T))\le \frac{C_0 \lambda_i^{-1}}{T - (T+\lambda_i^{-1}t)} = \frac{C_0}{-t}.$$

By the compactness Theorem (Theorem 7.5 \cite{L}) for coupled Ricci flows, combined with the $\kappa$-noncollapsed condition, there exists a subsequence
$$(M,g_i(t),\f_i(t),p) \to  (M_\infty,g_\infty(t),\f_\infty, p_\infty),$$
converging in the Cheeger-Gromov sense, where $(g_\infty(t), \f_\infty(t))$ still satisfies the coupled Ricci flow equation (\ref{RFHE}).

By the gradient estimate (\ref{List estimate 0}) of $\phi(t)$, 
$$\sup_M \abs{\nabla \f(t)}_{g(t)}\le \frac{\sup |\f(0)|^2}{t} = \frac{C^2}{t}.$$
In terms of the rescaled solutions $(g_i(t), \f_i(t))$, 
$$\sup_M \abs{\nabla \f_i(t)}_{g_i(t)} = \lambda_i^{-1}\sup_M \abs{\nabla \f(\lambda_i^{-1}t + T)}_{g(\lambda_i^{-1}t + T)}\le \lambda_i^{-1} \frac{C^2}{\lambda_i^{-1} t+T} = \frac{C^2}{t+\lambda_i T}.$$
By the $C^\infty$ convergence of $(g_i,\f_i)$ it follows that
$$\sup_{M_\infty}\abs{\nabla \f_\infty(t)} = 0,$$
hence $\f_\infty= $ const, and the coupled Ricci flow $(g_\infty,\f_\infty)$ becomes the standard Ricci flow.

From now on, we will denote the ``reduced'' function $\ell_\infty$ constructed in the previous section by $\ell_p: M\times [0,T)\to \R$. Denote 
$$\ell_p^i: (M\times [-\lambda_i T, 0), g_i, \f_i, p)\to \R$$
by $\ell^i_p (q,t) = \ell_p(q,\lambda_i^{-1} t + T)$. By Lemma \ref{lemma estimate}, $\ell^i_p$ satisfy the following estimates
$$C(n,C_0)^{-1} \frac{d^2_{g_i(t)}(p,q)}{-t} - C(n,C_0) \le \ell^i_p(q,t)\le  C(n,C_0) \frac{d^2_{g_i(t)}(p,q)}{-t} + C(n,C_0),$$
$$\abs{\nabla \ell^i_p}(q,t)\le \frac{C(n,C_0)}{-t}\bk{ \frac{d^2_{g_i(t)}(p,q)}{-t}  + 1 }$$
and 
$$\Big|\frac{\partial}{\partial t}\ell^i_p(q,t)  \Big| \le \frac{C(n,C_0)}{-t}\bk{ \frac{d^2_{g_i(t)}(p,q)}{-t}  + 1 },$$
for any $(q,t)\in M\times [-\lambda_i T, 0)$, and the last two inequalities are understood as the derivatives of locally Lipschitz functions. Note that the constant $C(n,C_0)$ above is independent of $i$.

So we can extract a subsequence of $\ell^i_p$ which converges in the Cheeger-Gromov-$C^0_{loc}(M_\infty\times (-\infty,0))$ sense to a function $\ell^\infty: M_\infty \times (-\infty, 0)\to \R$, and the convergence is uniform on any compact subset of $M_\infty\times (-\infty, 0)$. Define
\bea
V^\infty(t)=(4\pi(-t))^{-n/2} \int_{M_\infty} e^{-\ell^\infty(q,t)} dV_{g_\infty(t)}.
\eea

The reduced volume functions
\begin{equation}\label{reduce vol function}V^i(t) = (4\pi(-t))^{-n/2}\int_M e^{-\ell^i_p(q,t)  } dV_{g_i(t)}(q) = V_\infty(\lambda_i^{-1} t + T)\end{equation}
are nondecreasing for each fixed $t$ since $V_\infty(s)$ is nondecreasing in $s\in [0,T)$ and bounded above by $1$. So for any fixed $t<0$, we have
\begin{equation}\label{eqn converge}
\lim_{i\to \infty} V^i(t) = \lim_{t'\to T} V_\infty(t'),
\end{equation}
noting that the RHS is independent of $t$. 
On the other hand, by the dominated convergence theorem and smooth convergence of $g_i(t)$ to $g_\infty(t)$, we have
\begin{equation}
	\lim_{i\to\infty} V^i(t) = \lim_{i\to \infty}(4\pi (-t))^{-n/2} \int_M e^{-\ell^i_p(q,t)} dV_{g_i(t)}  = V^\infty(t),
\end{equation}
Hence $V^\infty(t)$ is constant and independent of $t$. By Lemma \ref{lemma 15}, it follows that the limit metric $g_\infty(t)$ is a gradient shrinking Ricci soliton.

\subsection{Non-triviality of the soliton}

It remains only to show the Ricci soliton $g_\infty$ is not flat, when the base point $p$ is a Type I singularity of the flow $(g(t), \phi(t))$. We will use Perelman's pseudo locality and a contradiction argument. First recall that the pseudolocality states that there exists a uniform constant $\varepsilon_0=\varepsilon_0(n)$ such that if the extended Ricci flow solution $(M,g(t),\phi(t))$ satisfies $\|\phi(0)\|_{L^\infty ( B_{g(0)}(p,1) )}\le C$ and the geometry of $( B_{g(0)}(p,1),g(0)  )$ is sufficiently close to that of the Euclidean unit ball $(B(0,1),g_{can})$ in the $C^2$ sense, then we have
$$|Rm|(x,g(t))\le 10\varepsilon_0^{-2} , \quad \forall t\in (\frac{\varepsilon_0^2}{2},\varepsilon_0^2),~\forall x\in B_{g(\varepsilon_0^2/2)}(p, \varepsilon_0).$$
Suppose the limit metric $g_\infty(t)$ is flat, and by assumption it is also $\kappa$-noncollapsed, so by the uniformization theorem $(M_\infty,g_\infty)\cong (\R^n,g_{can})$. By the smooth convergence of $(M,g_i(t),\phi_i(t),p)$ to $(M_\infty,g_\infty(t))$, it follows that when $i$ is large enough, the three conditions in the pseudolocality theorem are satisfied on $( B_{g_i(-\varepsilon_0^2)}(p,1),g_i(-\varepsilon_0^2)  )$. So  it follows that
$$|Rm|(x,g_i(-\varepsilon_0^2 + t))\le 10\varepsilon_0^{-2},\quad \forall t\in (\frac{\varepsilon_0^2}{2},\varepsilon_0^2), ~\forall x\in B_{g_i(-\varepsilon_0^2/2)}(p,\varepsilon_0).$$
Fix a large $i$, and in terms of the original flow, we have
\begin{equation}\label{final estimate 1}
|Rm|(x,g(T+\lambda_i^{-1}(t-\varepsilon_0^2)))\le 10 \lambda_i\varepsilon_0^{-2},
\end{equation}
for all $t\in (\frac{\varepsilon_0^2}{2},\varepsilon_0^2)$ and $x\in B_{g(T-\lambda_i^{-1}\varepsilon_0^2/2)}(p,\frac{\varepsilon_0}{\sqrt{\lambda_i}})$.

However, since $p$ is a Type I singularity point, which means there exists a sequence of space-time points $(p_\alpha,t_\alpha)$ such that $t_\alpha\to T$, $p_\alpha\to p$ and 
\begin{equation}\label{final estimate 2}
|Rm|(p_\alpha, g(t_\alpha))\ge \frac{c}{T-t_\alpha}, \quad \mbox{for some }c>0.
\end{equation} 
However, by (\ref{final estimate 1}) when $\alpha$ is large enough
$$|Rm|(p_\alpha, g(t_\alpha))\le 10 \lambda_i \varepsilon_0^{-2},$$
and  the RHS is uniformly bounded above, and this contradicts (\ref{final estimate 2}).

\section*{Appendix}
\setcounter{equation}{0}
\setcounter{section}{7}
In this appendix, we provide a proof of Theorem \ref{coupled Zhang}, which extends Zhang's theorem (\cite{Z}) to the case of extended Ricci flow. We may assume $t= -1$ in (\ref{Ricci soliton}) and (\ref{Ricci soliton 2}) and for notational convenience we suppress the subindex $\infty$ in $(M_\infty,g_\infty,\phi_\infty,f_\infty)$. We fix a point $p\in M$.

\begin{lemma}\label{lemma 1}
We have 
$$S + \abs{\nabla f} - f \equiv \mathrm{const}.$$
\end{lemma}
{\it Proof.}
For a function $f$, we will denote $f_i = \nabla_i f$ and $f_{ij} = \nabla_i\nabla_j f$, etc. We will also use a comma to denote the covariant derivative, e.g. $R_{ij,k}=\nabla_k R_{ij}$. Taking covariant derivative for (\ref{Ricci soliton}), we have
$$R_{ij,j} -  \phi_{ij}\phi_j -  \phi_i\Delta\phi + f_{ijj} = 0,$$
by the Ricci identity
$$f_{ijj} = f_{jij} = f_{jji} + f_j R_{ij} = (\Delta f)_i + R_{ij}f_j,$$
and the contracted second Bianchi identity $R_{ij,j}=\frac 1 2 R_{,i}$, we get
$$\frac 1 2 R_{,i} -  \phi_i \Delta \phi - \frac{1}{2}(\abs{\nabla\phi})_{i} + (\Delta f)_i + R_{ij}f_j = 0,$$ 
combining with the soliton equation (\ref{Ricci soliton}),  $R_{ij} = \frac 1 2 g_{ij} - f_{ij} +  \phi_i\phi_j$ and $\Delta\phi = \phi_i f_i$, we have
$$\bk{ \frac 1 2 R - \frac 1 2 \abs{\nabla \phi} + \Delta f + \frac 1 2 f -\frac 1 2 \abs{\nabla f}   }_i = 0$$
so we conclude that 
$$\frac 1 2 R - \frac 1 2 \abs{\nabla \phi} + \Delta f + \frac 1 2 f -\frac 1 2 \abs{\nabla f}$$
is a constant. Combining with $S+\Delta f= \frac n 2$, we finish the proof of the lemma. 

Without loss of generality, we
normalize $f$ such that
$$
S + \abs{\nabla f}  - f = 0.
$$
 By direct calculations we have
\begin{equation}\label{eqn:10}
\Delta S = \langle\nabla f,\nabla S\rangle + S - 2\abs{\sic}- 2 (\Delta \phi)^2.
\end{equation}
For completeness, we give the proof of (\ref{eqn:10}).

Noting that by Lemma \ref{lemma 1} and Bochner formula,
$$\Delta S = \Delta f - 2\abs{\nabla^2 f} - 2\innpro{\nabla \Delta f,\nabla f} - 2\ric(\nabla f,\nabla f),$$
and $$-\innpro{\nabla f,\nabla S}= - \abs{\nabla f} + \innpro{\nabla f, \nabla \abs{\nabla f}},$$
moreover, by the soliton equation (\ref{Ricci soliton}) we have
$$\abs{\nabla^2f} = \abs{\sic}+ \frac n 4 - S, $$
combining with the equations above, it follows that
$$ \Delta S - \innpro{\nabla f,\nabla S}= S - 2\abs{\sic} - 2 \ric(\nabla f,\nabla f) + \abs{\nabla f} - \innpro{\nabla f,\nabla \abs{\nabla f}},
$$
and the last three terms on RHS are equal to
$$-2(R_{ij}+f_{ij})f_if_j + \abs{\nabla f} = -2(\frac 1 2 g_{ij}+ \phi_i\phi_j)f_if_j + \abs{\nabla f}= -2 (\innpro{\nabla\phi,\nabla f})^2$$
now the formula (\ref{eqn:10}) follows from $\Delta \phi = \innpro{\nabla f,\nabla \phi}$.

Now we go back to the proof of nonnegativity of $S$ for the complete soliton (\ref{Ricci soliton}). We rewrite the estimate (\ref{Perelman estimate}) as follows: if $\ric\le (n-1)K$ in the ball $B(p,r_0)$, then
$$
\Delta d(p,\cdot)\le (n-1)(r_0^{-1} + \frac{2}{3}K r_0) -\int_\gamma \ric(\gamma',\gamma')ds,
$$
where $\gamma(s)$ is a normal geodesic initiating at $p$. Recall the equation $\sic + \nabla^2 f = 1/2 g$, $\sic(\gamma',\gamma') + \nabla^2f(\gamma',\gamma') = 1/2$, and 
$$
\frac{d}{ds} f(\gamma(s)) = \langle \nabla f, \gamma'\rangle,
$$
$$
\frac{d^2}{ds^2} f(\gamma(s)) = \langle \nabla_{\gamma'}\nabla f,\gamma'\rangle = \nabla^2 f(\gamma',\gamma'),
$$
thus 
\bea
\int_\gamma \ric(\gamma',\gamma')ds& = &\frac{1}{2} d(p,x) - \int_0^{d(p,x)} \frac{d^2}{ds^2} f(\gamma(s)) ds+\int_0^{d(p,x)} d\phi\otimes d\phi (\gamma',\gamma')ds \nonumber \\
&\ge& \frac{1}{2}d(p,x) - \langle \nabla f(\gamma(s)),\gamma'(s)\rangle |_{s=0}^{s=d(p,x)} \nonumber\\
&\ge& \frac{1}{2}d(p,x) - \langle\nabla f(x), \nabla_x d(p,x)\rangle - |\nabla f(p)|.\nonumber
\eea
Hence, we have for $d(x) = d(p,x)$
$$
\Delta d - \langle \nabla f,\nabla d\rangle\le (n-1)(\frac{2}{3}Kr_0+ r_0^{-1}) - \frac{1}{2} d + |\nabla f|(p).
$$
For any fixed $p\in M$, note that $(M,g)$ is complete noncompact. There exists an $r_0=r_0(p)>0$, such that $\ric(g)\le (n-1)r_0^{-2}$ on the ball $B(p,r_0)$. Choose a cut-off function $\eta(x)$ such that it is $1$ when $x$ is less than equal to $1$, vanishes when $x$ is greater than $2$, and $|\eta'|^2/\eta\le 4$, $\eta'\le 0$ and $|\eta''|, |\eta'|\le 2$ in ${\mathbf R}$. Then define a new function on $M$ by (by abusing notation we still use $\eta$ to denote the new function)
$$
\eta(x):= \eta\bk{\frac{d(p,x)}{Ar_0}},
$$
where $A>0$ is a large constant. Define $u=S\eta$, and we calculate
\bea 
\Delta u =  \eta \Delta S + \frac{S\eta'}{Ar_0}\Delta d + \frac{S\eta''}{(Ar_0)^2} + 2\langle \nabla S,\nabla \eta\rangle, 
\eea
We consider two cases $u_{\min}\ge 0$ or $u_{\min}<0$. For the first case, we are already done, since that means $S(p) = u(p)\ge 0$; in the latter case, we know the minimum must be inside the support of $\eta$, which is a point in $B(p,Ar_0)$ denoted by $p_{\min}$. At the minimum point $p_{\min}$, we have $\nabla u = 0$ and $\Delta u\ge 0$, hence at $p_{\min}$, we see $\nabla S =\frac{-S\nabla\eta}{\eta}$ and $u\eta'\ge 0$, and
\bea \Delta u & = &\eta\Big(\langle \nabla S,\nabla f\rangle + S - 2|\sic|^2 - 2(\Delta\phi)^2\Big) + \frac{u \eta'}{Ar_0\eta} \Delta d + \frac{u\eta''}{(Ar_0)^2\eta} + 2\langle \nabla S,\nabla \eta\rangle \nonumber\\
& \le& \eta\Big(\langle \nabla S,\nabla f\rangle + S - 2|\sic|^2\Big) + \frac{u \eta'}{Ar_0\eta} \Delta d + \frac{u\eta''}{(Ar_0)^2\eta} - 2S\frac{|\nabla \eta|^2}{\eta}\nonumber\\
&=&\eta\Big(\langle \nabla S,\nabla f\rangle + S - 2|\sic|^2\Big) + \frac{u \eta'}{Ar_0\eta} \Delta d + \frac{u\eta''}{(Ar_0)^2\eta} - 2u \frac{|\eta'|^2}{\eta^2}\frac{1}{(Ar_0)^2}\nonumber\\
&\le& \eta\Big(\langle \nabla S,\nabla f\rangle + S - 2|\sic|^2\Big) + \frac{u\eta'}{Ar_0\eta} (\langle\nabla f,\nabla d\rangle + 2nr_0^{-1} + |\nabla f|(p)) -2u \frac{|\eta'|^2}{\eta^2}\frac{1}{(Ar_0)^2}.\nonumber
\eea
Observing that at $p_{\min}$, the first term in RHS is $-\eta S \frac{\innpro{\nabla \eta,\nabla f}}{\eta} = - \frac{u\eta'}{Ar_0\eta}\innpro{\nabla f,\nabla d}$, which is cancelled with one term in RHS.
Hence we have at the minimum point $p_{\min}$, 
\begin{equation}\label{eqn:1}
 u - 2\eta |\sic|^2 + \frac{u\eta'}{Ar_0^2\eta} + \frac{2nu\eta'}{Ar_0\eta}|\nabla f|(p) - 2u\frac{|\eta'|^2}{\eta^2}\frac{1}{(Ar_0)^2}\ge 0.
\end{equation}
On the other hand,  by Cauchy-Schwarz $-2\eta|\sic|^2\le -2\frac{\eta S^2}{n} = -\frac{2u^2}{n\eta}$, multiplying by $\eta$ in (\ref{eqn:1}), we see that
$$
\eta u - \frac{2u^2}{n} + \frac{2nu\eta'}{Ar_0^2} + \frac{u\eta'}{Ar_0}|\nabla f|(p) - 2u\frac{|\eta'|^2/\eta}{(Ar_0)^2}\ge 0,
$$
that is
$$
0\le -\frac{2u^2}{n} - |u| + \frac{4n|u|}{Ar_0^2} + \frac{2|u|}{Ar_0}|\nabla f|(p) + 4\frac{|u|}{(Ar_0)^2},
$$
thus we conclude that at the minimum point $p_{\min}$,
$$
|u|\le \frac{C(|\nabla f|(p),n)}{A r_0^2},
$$
for any constant $A>0$, and this implies $$S(p)=u(p)\ge u_{\min}\ge - \frac{C(n,|\nabla f|(p))}{Ar_0^2}$$ Note that in this estimate, the constants $C(n,|\nabla f(p)|)$ and $r_0$ only depend on the chosen point $p$ and $|\nabla f(p)|$, and are independent of $A$, so letting $A\to \infty$, we have $S(p) = u(p)\ge 0$. Since $p$ is an arbitrary point in $M$, we conclude that $S\ge 0$ on $M$.

From $S+ \abs{\nabla f} = f$, we conclude that $ |\nabla \sqrt{f}|\le 2$, thus $|\nabla f|\le C d(p,x)$ outside a compact set. Hence, the completeness of the metric $g$ implies the completeness of the vector field $\nabla f$, which means it generates a one parameter family of diffeomorphisms.

\bigskip

\noindent{\bf Acknowledgement:} Bin Guo would like to thank Prof. Jian Song for his support and encouragement, and many helpful discussions. The authors would like to thank Prof. Reto M\"uller for bringing some additional references to their attention, and to Prof. Peter Topping for clarifying for them some differences between the results in \cite{EMT} and in \cite{CZ}.

\bigskip

\noindent Department of Mathematics, Columbia University, New York, NY 10027, USA\\
bguo@math.columbia.edu, zjhuang@math.columbia.edu, phong@math.columbia.edu

\end{document}